\documentclass[12pt]{article}
\usepackage{amsmath}
\usepackage{amssymb}
\usepackage{latexsym}
\usepackage{amsfonts}
\usepackage{graphicx}
\usepackage{color}
\usepackage{tikz}
\usepackage{lineno}

\newcommand{\gtr}[1]{\gamma_{[3R]}\left( #1 \right)}
\newtheorem{thm}{Theorem}
\newtheorem{obs}[thm]{Observation}

\newtheorem{cor}[thm]{Corollary}
\newtheorem{prop}[thm]{Proposition}

\DeclareMathOperator{\opt}{\emph{Opt}}

\begin{document}

\title{Triple Roman Domination in Graphs}
\date{}
\author{$^{(1)}$H. Abdollahzadeh Ahangar, $^{(2)}$M.P. \'Alvarez, $^{(3)}$M.
Chellali, \\
$^{(4)}$S.M. Sheikholeslami and $^{(5)}$J.C. Valenzuela-Tripodoro\vspace{%
7.5mm} \\
$^{(1)}${\small Department of Mathematics}\\
{\small Babol Noshirvani University of Technology}\\
{\small Shariati Ave., Babol, I.R. Iran, Postal Code: 47148-71167}\\
{\small ha.ahangar@nit.ac.ir\vspace{2mm}}\\
$^{(2)}${\small Department of Statistic and O.R., University of C\'adiz,
Spain.}\\
{\small pilar.ruiz@uca.es \vspace{2mm}}\\
$^{(3)}${\small LAMDA-RO Laboratory, Department of Mathematics}\\
{\small University of Blida}\\
{\small Blida, Algeria}\\
{\small m\_chellali@yahoo.com\vspace{2mm}}\\
$^{(4)}${\small Department of Mathematics}\\
{\small Azarbaijan Shahid Madani University}\\
{\small Tabriz, Iran}\\
{\small s.m.sheikholeslami@azaruniv.ac.ir\vspace{2mm} } \\
$^{(5)}${\small Department of Mathematics, University of C\'adiz, Spain.} \\
{\small jcarlos.valenzuela@uca.es\vspace{5mm} }}
\maketitle

\begin{abstract}
The Roman domination in graphs is well-studied in graph theory. The topic is
related to a defensive strategy problem in which the Roman legions are
settled in some secure cities of the Roman Empire. The deployment of the
legions around the Empire is designed in such a way that a sudden attack to
any undefended city could be quelled by a legion from a strong neighbour.
There is an additional condition: no legion can move if doing so leaves its
base city defenceless. In this {manuscript} we start the study of a variant
of Roman domination in graphs: the triple Roman domination. We consider that
any city of the Roman Empire must be able to be defended by {at least} three
legions. These legions should be either in the attacked city or in one of
its neighbours. {We determine various bounds on the triple Roman domination
number for general graphs, and we give exact values for some graph families.
Moreover, complexity results are also obtained. }

\end{abstract}

\section{Introduction}

Only simple, undirected and {non-trivial} connected graphs {(ntc-graphs)}
will be considered in this {manuscript}. The set of vertices of the graph ${%
\Gamma }$ is denoted by $V=V({\Gamma })$ and the edge set is $E=E({\Gamma })$%
. The order of a graph is the number of vertices of the graph $\Gamma $ and
it is denoted by $p=p({\Gamma )}$. The size of $\Gamma $ is the cardinality
of the edge set and it is denoted by $q=q({\Gamma )}.$ An edge joining the
vertices $u$ and $v$ is denoted by $e=uv$ and it is said to be incident with
both end-vertices. For every vertex $v\in V$, the number of edges that are
incident with $v$ is the degree of $v$ and is denoted by $d_{\Gamma }(v)$,
or simply by $d(v)$ if the underlying graph is clear. Given two sets of
vertices $S,T\subseteq V,$ we denote by $E[S,T]$ the set of edges having one
end-vertex in $S$ and the other end-vertex in $T.$ The \emph{open {neighbour}%
hood} $N(v)$ is the set $\{u\in V(\Gamma ):uv\in E(\Gamma )\}$ and the \emph{%
closed {neighbour}hood} of $v$ is the set $N[v]=N(v)\cup \{v\}$. The \emph{%
minimum} and \emph{maximum degree} of a graph $\Gamma $ are denoted by $%
\delta =\delta (\Gamma )$ and $\Delta =\Delta (\Gamma )$, respectively. Any
vertex of degree one is called a \textit{leaf }and a \textit{{stem} vertex}
is a vertex {adjacent to }a leaf. {For a stem }$v$, let $L_{v}$ {be the set
of leaves adjacent to }$v.$

{A }\textit{path} {$P_{t+1}$} in a graph is a set of different vertices $%
\{v_{i}:i=0,\ldots ,t\}$ and the set of edges $\{v_{i}v_{i+1}:i=0,\ldots
,t-1\}$. A \textit{cycle} {$C_{t}$} is a path in which $v_{0}=v_{t}$. A
\textit{connected graph} is a graph in which any pair or vertices could be
joined by a path. {The \textit{distance} $d_{G}(u,v)$ between two vertices $u
$ and $v$ in a connected graph $G$ is the length of a shortest $u-v$ path in
$G.$ The \emph{diameter} of $G$, denoted by $\mathrm{diam}(G)$, is the
maximum value among distances between all pair of vertices of $G$. The \emph{%
girth} of $G$, denoted by $g(G)$, is the minimum length of a cycle in $G$. }

{A \emph{tree} is an acyclic connected graph. A tree }$T$ is a \emph{double
star} {if it contains exactly two vertices that are not leaves}.\ A double
star with, respectively $r$ and $s$ leaves attached at each stem is denoted $%
DS_{r,s}.$ {For a vertex $v$ in a rooted tree $T$, let $C(v)$ and $D(v)$
denote the set of children and descendants of $v$, respectively and let $%
D[v]=D(v)\cup \{v\}$. Also, the depth of $v$, $
\mathrm{depth}(v)$, is the largest distance from $v$ to a vertex in $D(v)$.
The \emph{maximal subtree} \emph{rooted at $v$}, denoted by $T_{v}$,
consists of $v$ and all its descendants.}

A \textit{dominating set }{(DS)} in a graph $\Gamma $ is a set of vertices $%
S\subseteq V(\Gamma )$ such that any vertex of $V-S$ is {connecting} to {at
least} one vertex of $S.$ The \textit{domination number} $\gamma (\Gamma )$
equals the minimum cardinality of a {DS} in $\Gamma $.

Let $h:V(\Gamma)\rightarrow \{0,1,2,\ldots,k\}$ be a function, and let $%
(V_{0},V_{1}, V_{2},\ldots,V_{k})$ be the ordered partition of $V=V(\Gamma)$
induced by $h$, where $V_{i}=\{v\in V:h(v)=i\}$ for $i\in\{0,1,\ldots,k\}$.
There is a $1$-$1$ correspondence between the functions $h:V\rightarrow%
\{0,1,2,\ldots,k\}$ and the ordered partitions $(V_{0},V_{1},V_{2},%
\ldots,V_{k}) $ of $V$, so we will write $h=(V_{0},V_{1},V_{2},\ldots,V_{k})$%
. {\ For any subset $A\subseteq V$ we denote by $h_{|A}:A\rightarrow
\{0,1,\ldots, k \}$ the restriction of the function $h$ to the smaller
domain $A$.}

A function $h:V(\Gamma)\rightarrow \{0,1,2\}$ is a \textit{Roman dominating
function} (RDF) on $\Gamma$ if every vertex $u\in V$ for which $h(u)=0$ is {%
connecting to} {at least} one vertex $v$ for which $h(v)=2$. The weight of {a%
} RDF is the value $h(V(\Gamma))=\sum_{u\in V(\Gamma)}h(u).$ The \textit{%
Roman domination number} {RD-number} $\gamma _{R}(\Gamma)$ is the minimum
weight of an RDF on $\Gamma$. The concept of Roman domination was introduced
by Cockayne et al. in \cite{CDHH} and was inspired by the {manuscript} of {%
ReVelle and Rosing} \cite{RR}, and Stewart~\cite{S} about the
defensive strategy of the Roman Empire decreed by Constantine I The Great.
The strategy establishes that any undefended place should {be a neighbour}
to a 'strong' city in which there are (exactly) 2 legions deployed. In that
way, one of these two legions could quell any sudden attack to the
undefended city, without leaving its base site undefended. {Clearly}, taking
into account the defensive strategy decreed by the Roman Emperor, the aim is
to minimize the total cost of the legions deployed along the empire.

The topic of {RD-number} in graphs {has been} extensively studied for the
last 15 years, where many {manuscript}s have been published. Moreover,
several new variations of {RD-number} were introduced. For example, weak
Roman domination \cite{HH}, independent Roman domination \cite{AEJM},
maximal Roman domination \cite{ASSTV}, mixed Roman domination \cite{AHV},
Roman $\{2\}$-domination \cite{CHHA} and recently total Roman $\{2\}$%
-domination \cite{acsv}, for more see \cite{acs1, acs2, HW}.

In 2016, {Beeler et al.} \cite{BHH} introduced the double Roman
domination number {DRD-number}, by proposing a stronger version of the {%
RD-number} in which any city of the empire could be defended by {at least}
two legions. In order to achieve this goal, up to three legions could be
placed in a given location. This provides a defensive capacity which doubles
the one of the original strategy and, simultaneously, does not increase the
total cost as much as one would expect. Formally, a \textit{double Roman
dominating function} {DRDF} in a graph $\Gamma$ is a function $%
h:V(\Gamma)\rightarrow \{0,1,2,3\}$ that satisfies the following conditions:
\vspace{-0.25cm}

\begin{enumerate}
\item[a)] If $h(v)=0,$ then vertex $v$ must have {at least} two neighbours
in $V_{2}$ or {at least} one neighbour in $V_{3}.$

\item[b)] If $h(v)=1,$ then vertex $v$ must have {at least} one {neighbour}
in $V_{2}\cup V_{3}.$
\end{enumerate}

The minimum weight $h(V(\Gamma))=\sum_{v\in V}h(v)$ of a {DRDF} is the
\textit{DRD-number} $\gamma _{dR}(\Gamma)$ of the graph $\Gamma.$ For more
see {\cite{aaacs, AA, acs, cfsy,v2, K, z,v1}}.

Beeler et al. pointed out in \cite{BHH} that sometimes doubling the defence
capacity (that is to say, defend against {any} attack with {at least} two
legions), only increases the cost up to, {at most}, 50\%, as can be easily
verified for the star $K_{1,p-1}$ and the complete bipartite graph $%
K_{2,p-2}.$

In this {manuscript}, we introduce a generalization of the {DRD-number} in
which we assume that any undefended place could be defended from a sudden
attack with, {at least}, $k$ legions without leaving any '{neighbour}ing
strong-city' without military forces. We can think of this generalization as
in {an \emph{umpteenth Roman Domination}} or a \emph{k-th Roman Domination}.

More precisely, let $h$ be a function that assigns labels from the set $%
\{0,1,\ldots ,k+1\}$ to the vertices of the graph $\Gamma$. Given a vertex $%
v\in V(\Gamma),$ the \emph{active {neighbour}hood} of $v$, denoted by $AN(v)$%
, is the set of vertices $w\in N_{\Gamma}(v)$ such that $h(w)\geq 1.$ Let $%
AN[v]=AN(v)\cup \{v\}.$

A [k]-\textit{RDF} is a function $h:V\rightarrow \{0,1,\ldots ,k+1\}$ such
that for {any} vertex $v\in V$ with $h(v)<k$,
\begin{equation*}
h\left( AN[v]\right) \geq |AN(v)|+k.
\end{equation*}%
The weight of a [k]-{RDF} is the value $h(V)=\sum_{v\in V}h(v)$, and the
minimum weight of such a function [k]-\textit{RD-number} of $\Gamma$,
denoted by $\gamma _{\lbrack kR]}(\Gamma)$.

Let us point out that for $k=2$ the previous definition matches that of the {%
DRD-number}. In this {manuscript} we focus our attention to the \emph{triple
Roman domination number} {TRD-number} case, so that for {any} vertex $v\in V$
with $h(v)<3 $, it must happen that
\begin{equation*}
h\left( AN[v]\right) \geq |AN(v)|+3.
\end{equation*}%
That is to say, we have to label the vertices of the graph, with labels from
the set $\{0,1,2,3,4\}$, so that:

\begin{enumerate}
\item If $h(v)=0,$ then $v$ must have either one neighbour in $V_{4}$, or
either two neighbours in $V_{2}\cup V_{3}$ (one neighbour in $V_{2}$ and
another one in $V_{3}$) or either three neighbours in $V_{2}.$

\item If $h(v)=1,$ then $v$ must have either one neighbour in $V_{3}\cup
V_{4}$ or either two neighbours in $V_{2}.$

\item If $h(v)=2,$ then $v$ must have one neighbour in $V_{2}\cup V_{3}\cup
V_{4}.$
\end{enumerate}

Let us denote by $\gamma _{\lbrack 3R]}\left( \Gamma\right) $ the minimum
weight of a triple Roman domination function (3RDF) in $\Gamma$.

\begin{figure}[h]
\centering
\begin{tikzpicture}[scale=.85, transform shape]

\draw (0,3) ellipse (0.7cm and 1.75cm) ;
\node [draw, shape=circle,fill=black,scale=0.5] (u1) at  (0,2) {};
\node [draw, shape=circle,fill=black,scale=0.5] (u2) at  (0,3) {};
\node [draw, shape=circle,fill=black,scale=0.5] (u3) at  (0,4) {};
\draw (2,3) ellipse (0.75cm and 2cm);
\node [draw, shape=circle,fill=black,scale=0.5] (v1) at  (2,1.5) {};
\node [draw, shape=circle,fill=black,scale=0.5] (v2) at  (2,2.5) {};
\node [draw, shape=circle,fill=black,scale=0.5] (v3) at  (2,3.5) {};
\node [draw, shape=circle,fill=black,scale=0.5] (v4) at  (2,4.5) {};
\draw (u1)--(v1);
\draw (u1)--(v2);
\draw (u1)--(v3);
\draw (u1)--(v4);
\draw (u2)--(v1);
\draw (u2)--(v2);
\draw (u2)--(v3);
\draw (u2)--(v4);
\draw (u3)--(v1);
\draw (u3)--(v2);
\draw (u3)--(v3);
\draw (u3)--(v4);
\node[draw=none, fill=none] at (1,.5){ $\gamma_R(K_{3,4})=4 $};
\node[draw=none, fill=none] at (-0.25,4){ $2$};
\node[draw=none, fill=none] at (-0.25,3){ $0$};
\node[draw=none, fill=none] at (-0.25,2){ $0$};
\node[draw=none, fill=none] at (2.25,4.5){ $2$};
\node[draw=none, fill=none] at (2.25,3.5){ $0$};
\node[draw=none, fill=none] at (2.25,2.5){ $0$};
\node[draw=none, fill=none] at (2.25,1.5){ $0$};
\draw (4,3) ellipse (0.7cm and 1.75cm) ;
\node [draw, shape=circle,fill=black,scale=0.5] (w1) at  (4,2) {};
\node [draw, shape=circle,fill=black,scale=0.5] (w2) at  (4,3) {};
\node [draw, shape=circle,fill=black,scale=0.5] (w3) at  (4,4) {};
\draw (6,3) ellipse (0.75cm and 2cm) ;
\node [draw, shape=circle,fill=black,scale=0.5] (w4) at  (6,1.5) {};
\node [draw, shape=circle,fill=black,scale=0.5] (w5) at  (6,2.5) {};
\node [draw, shape=circle,fill=black,scale=0.5] (w6) at  (6,3.5) {};
\node [draw, shape=circle,fill=black,scale=0.5] (w7) at  (6,4.5) {};
\draw (w1)--(w4);
\draw (w1)--(w5);
\draw (w1)--(w6);
\draw (w1)--(w7);
\draw (w2)--(w4);
\draw (w2)--(w5);
\draw (w2)--(w6);
\draw (w2)--(w7);
\draw (w3)--(w4);
\draw (w3)--(w5);
\draw (w3)--(w6);
\draw (w3)--(w7);
\node[draw=none, fill=none] at (5,.5){ $\gtr {K_{3,4}}=8 $};
\node[draw=none, fill=none] at (3.75,4){ $4$};
\node[draw=none, fill=none] at (3.75,3){ $0$};
\node[draw=none, fill=none] at (3.75,2){ $0$};
\node[draw=none, fill=none] at (6.25,4.5){ $4$};
\node[draw=none, fill=none] at (6.25,3.5){ $0$};
\node[draw=none, fill=none] at (6.25,2.5){ $0$};
\node[draw=none, fill=none] at (6.25,1.5){ $0$};
\draw (7.5,0)--(7.5,6);
\node [draw, shape=circle, fill=black, scale=0.5] (a1) at  (9.75 ,5 ) {};
\node [draw, shape=circle, fill=black, scale=0.5] (a2) at  (8 ,3.7) {};
\node [draw, shape=circle, fill=black, scale=0.5] (a3) at  (11.5 ,3.7) {};
\node [draw, shape=circle, fill=black, scale=0.5] (a4) at  (8.75 ,1.5 ) {};
\node [draw, shape=circle, fill=black, scale=0.5] (a5) at  (10.75 ,1.5 ) {};
\node [draw=none, fill=none] at  (9.75 , 5.45 ) {2};
\node [draw=none, fill=none] at  (8 ,4.15) {0};
\node [draw=none, fill=none] at  (11.5 ,4.15 ) {0};
\node [draw=none, fill=none] at  (8.4 ,1.8 ) {1};
\node [draw=none, fill=none] at  (11.1,1.8) {1};
\draw (a1)--(a2)--(a4)--(a5)--(a3)--(a1);
\node[draw=none, fill=none] at (9.75,.5){ $\gamma_{R}(C_{5})=4 $};
\node [draw, shape=circle, fill=black, scale=0.5] (b1) at  (14 ,5 ) {};
\node [draw, shape=circle, fill=black, scale=0.5] (b2) at  (12.25 ,3.7) {};
\node [draw, shape=circle, fill=black, scale=0.5] (b3) at  (15.75 ,3.7 ) {};
\node [draw, shape=circle, fill=black, scale=0.5] (b4) at  (13 ,1.5 ) {};
\node [draw, shape=circle, fill=black, scale=0.5] (b5) at  (15 ,1.5 ) {};
\node [draw=none, fill=none] at  (14 , 5.45 ) {3};
\node [draw=none, fill=none] at  (12.25 ,4.15) {0};
\node [draw=none, fill=none] at  (15.75 ,4.15 ) {0};
\node [draw=none, fill=none] at  (12.7 ,1.8 ) {2};
\node [draw=none, fill=none] at  (15.3 ,1.8) {2};
\draw (b1)--(b2)--(b4)--(b5)--(b3)--(b1);
\node[draw=none, fill=none] at (14,.5){ $\gtr{C_{5}}=7 $};

\end{tikzpicture}
\caption{{\protect\small Increasing the defence up to three times does not
always increase the cost to triple.}}
\label{bipypent}
\end{figure}
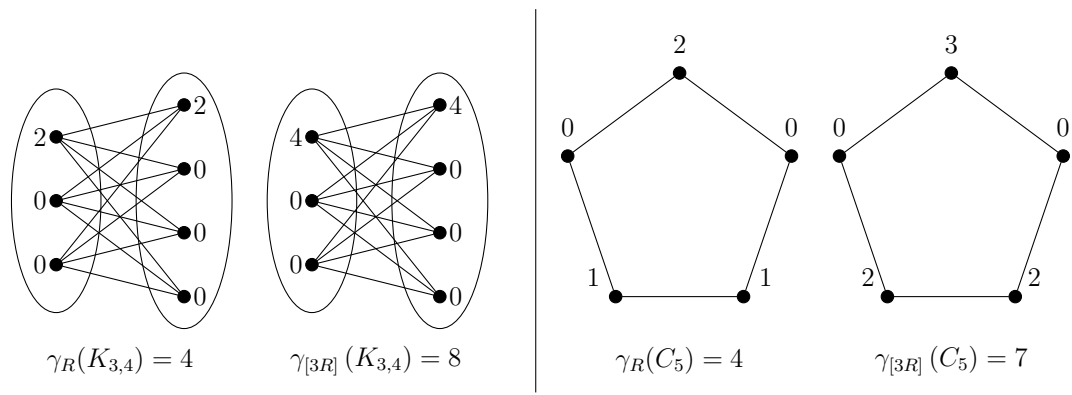

It is worth noting that there is a significant benefit in considering the {%
TRD-number} strategy in some cases. For example, the {RD-number} for the
bipartite graph $K_{3,4}$ equals $4$ while the corresponding {TRD-number}
would be only $8$ (see Fig. \ref{bipypent}). That is to say, we are able to
increase the defensive capacity up to three times whilst the total cost is
only double.

Moreover, in the case of the pentagon (see the right side in Fig. \ref%
{bipypent}), we need $4$ legions to defend the entire graph under the {%
RD-number} ($\gamma_R(C_5)=4$) and we only need $7$ legions to defend any
place, with {at least} three legions, according to the definition of the {%
TRD-number} ($\gamma_{[3R]}\left( C_5 \right)=7$). This means that we may
increase the defence to 300\%, because any sudden attack to an unsafe place
could be defended with {at least} three legions, with only an increase in
the cost of 75\%, since the number of legions needed goes from $4 $ to $7$.

In this {manuscript}, we initiate the study of {TRD-number}. We first show
that the problem of computing $\gamma_{\lbrack 3R]}(\Gamma)$ is NP-complete
for bipartite and chordal graphs. {Moreover, we show that it is possible to
compute this parameter in linear time for bounded clique-width graphs
including the class of trees. Also, we }establish various bounds on the {%
TRD-number} for general graphs. In particular, we show that for any {%
ntc-graph} $\Gamma$ of order $p\geq 2$, $\gamma _{\lbrack 3R]}(\Gamma)\leq
\frac{7p}{4}$, and we characterize the {ntc-graphs} attaining the upper
bound. Finally the exact values of the {TRD-number} for some graph families
are given.%


\section{Complexity results}

Our {goal} in {here} is to establish the NP-complete result for the {%
TRD-number} problem in bipartite and chordal graphs.\newline
\newline
TRIPLE\ ROM-DOM\newline
\textbf{Instance}: Graph $\Gamma =(V,E)$, positive integer $k\leq |V|$.%
\newline
\textbf{Question}: Does $\Gamma $ have a {3RDF} of weight {at most} $k$?

\bigskip

We show that this problem is NP-complete by reducing the well-known
NP-complete problem, Exact-$3$-Cover (X3C), to TRIPLE ROM-DOM.\newline
\newline
EXACT $3$-Cover (X3C)\newline
\textbf{Instance}: A finite set $X$ with $|X|=3q$ and a collection $C$ of $3$%
-{member} subsets of $X$.\newline
\textbf{Question}: Is there a subcollection $C^{\prime}$ of $C$ {for which} {%
any member} of $X$ appears in {perfectly} one {member} of $C^{\prime}$?%
\newline

\begin{thm}
\label{complexity_bipartite}Problem TRIPLE ROM-DOM is NP-Complete for
bipartite graphs.
\end{thm}

\noindent \textbf{Proof.} {Clearly} TRIPLE ROM-DOM is a member of $\mathcal{%
NP}$, since we can check in polynomial time that a function $l:V\rightarrow
\{0,1,2,3,4\}$ has weight {at most} $k$ and is a {3RDF}. Given an instance $%
(X,C)$ of X3C with $X=\{x_{1},x_{2},\ldots ,x_{3q}\}$ and $%
C=\{C_{1},C_{2},\ldots ,C_{t}\}$.

\begin{figure}[h]
\centering
\begin{tikzpicture}
\node [draw, shape=circle,fill=black,scale=0.5] (y1) at  (0,1) {};
\node at (-0.4,1) {$y_1$};
\node [draw, shape=circle,fill=black,scale=0.5] (y2) at  (1.5,1) {};
\node at (1.1,1) {$y_2$};
\node [draw, shape=circle,fill=black,scale=0.5] (y3) at  (3,1) {};
\node at (2.6,1) {$y_3$};
\node [draw, shape=circle,fill=black,scale=0.5] (y4) at  (4.5,1) {};
\node at (4.1,1) {$y_4$};
\node [draw, shape=circle,fill=black,scale=0.5] (y5) at  (6,1) {};
\node at (5.6,1) {$y_5$};
\node [draw, shape=circle,fill=black,scale=0.5] (y6) at  (7.5,1) {};
\node at (7.1,1) {$y_6$};
\node [draw, shape=circle,fill=black,scale=0.5] (x1) at  (0,0) {};
\node at (-0.4,0.05) {$x_1$};
\node [draw, shape=circle,fill=black,scale=0.5] (x2) at  (1.5,0) {};
\node at (1.1,0.05) {$x_2$};
\node [draw, shape=circle,fill=black,scale=0.5] (x3) at  (3,0) {};
\node at (2.6,0.05) {$x_3$};
\node [draw, shape=circle,fill=black,scale=0.5] (x4) at  (4.5,0) {};
\node at (4.1,0.05) {$x_4$};
\node [draw, shape=circle,fill=black,scale=0.5] (x5) at  (6,0) {};
\node at (5.6,0.05) {$x_5$};
\node [draw, shape=circle,fill=black,scale=0.5] (x6) at  (7.5,0) {};
\node at (7.1,0.05) {$x_6$};
\node [draw, shape=circle,fill=black,scale=0.5] (c1) at  (.75,-2) {};
\node at (0.4,-2.05) {$c_1$};
\node [draw, shape=circle,fill=black,scale=0.5] (c2) at  (2.75,-2) {};
\node at (2.4,-2.05) {$c_2$};
\node [draw, shape=circle,fill=black,scale=0.5] (c3) at  (4.75,-2) {};
\node at (4.4,-2.05) {$c_3$};
\node [draw, shape=circle,fill=black,scale=0.5] (c4) at  (6.75,-2) {};
\node at (6.4,-2.05) {$c_4$};
\node [draw, shape=circle,fill=black,scale=0.5] (a1) at  (.75,-3) {};
\node [draw, shape=circle,fill=black,scale=0.5] (a11) at  (.75,-3.6) {};
\node [draw, shape=circle,fill=black,scale=0.5] (a12) at  (.45,-3.6) {};
\node [draw, shape=circle,fill=black,scale=0.5] (a13) at  (1.05,-3.6) {};
\draw(a13)--(a1)--(a12);
\node [draw, shape=circle,fill=black,scale=0.5] (a2) at  (2.75,-3) {};
\node [draw, shape=circle,fill=black,scale=0.5] (a21) at  (2.75,-3.6) {};
\node [draw, shape=circle,fill=black,scale=0.5] (a22) at  (2.45,-3.6) {};
\node [draw, shape=circle,fill=black,scale=0.5] (a23) at  (3.05,-3.6) {};
\draw(a23)--(a2)--(a22);
\node [draw, shape=circle,fill=black,scale=0.5] (a3) at  (4.75,-3) {};
\node [draw, shape=circle,fill=black,scale=0.5] (a31) at  (4.75,-3.6) {};
\node [draw, shape=circle,fill=black,scale=0.5] (a32) at  (4.45,-3.6) {};
\node [draw, shape=circle,fill=black,scale=0.5] (a33) at  (5.05,-3.6) {};
\draw(a33)--(a3)--(a32);
\node [draw, shape=circle,fill=black,scale=0.5] (a4) at  (6.75,-3) {};
\node [draw, shape=circle,fill=black,scale=0.5] (a41) at  (6.75,-3.6) {};
\node [draw, shape=circle,fill=black,scale=0.5] (a42) at  (6.45,-3.6) {};
\node [draw, shape=circle,fill=black,scale=0.5] (a43) at  (7.05,-3.6) {};
\draw(a43)--(a4)--(a42);
\draw(a11)--(c1)--(x1)--(y1)--(x1)--(c2)--(a21)--(c2)--(x5)--(y5)--(x5)--(c1)--(x4)--(y4)--(x4)--(c4)--(x3)--(y3)--(x3)--(c3)--(x2)--(y2)--(x2)--(c4)--(a41);
\draw(c2)--(x6)--(y6)--(x6)--(c3)--(a31);
\end{tikzpicture}
\caption{NP-completeness of 3RDF for bipartite.}
\label{fig-1}
\end{figure}
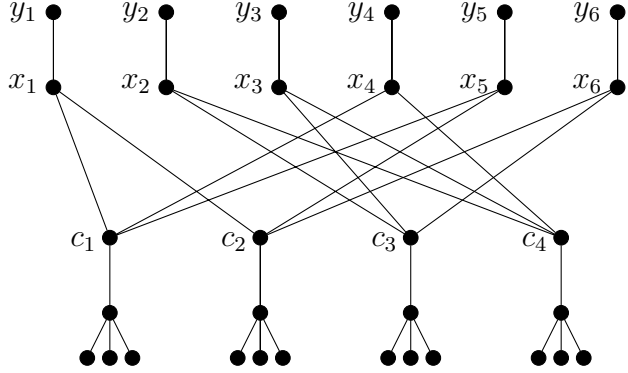

{\ We construct the bipartite graph $\Gamma$ as follows: }for {any} $%
x_{i}\in X$, we create a path $P_{2}:$ $x_{i}y_{i};$ for {any} $C_{j}\in C$
we build a star $H_{j}=K_{1,4}$ centered at $w_{j}$ with one of its leaves
labeled $c_{j}.$ Let $Z=\{c_{1},c_{2},\ldots,c_{t}\}.$ To achieve the
construction of $\Gamma,$ we add edges $c_{j}x_{i}$ {when} $x_{i}\in C_{j}$.
Set $k=4t+11q$.

{Assume} that the instance $X,C$ of X3C has a solution $C^{\prime }$. We
construct a {3RDF} $l$ on $\Gamma $ of weight $k$. We {label} a $4$ to all $%
w_{j}$'s, a $3$ to all $y_{i}$'s, a $0$ to all $x_{i}$'s and a $0 $ to
leaves of each $w_{j}.$ For every $c_{j},$ {label} a $2$ if $C_{j}\in
C^{\prime }$, and a $0$ if $C_{j}\notin C^{\prime }.$ Note that since $%
C^{\prime }$ exists, its cardinality is {exactly} $q$, and so the number of $%
c_{j}$'s with weight $2$ is $q$, having disjoint {neighbour}hoods in $%
\{x_{1},x_{2},...,x_{3q}\}$. Note that $C^{\prime }$ is a solution for X3C,
every vertex $x_{i}\in X$ satisfies $l(AN[x_{i}])\geq \left\vert
AN(x_{i})\right\vert +3.$ {Thus}, it is straightforward to {check} that $l$
is a {3RDF} with weight $l(V)=4t+2q+9q=k$.

Conversely, {let $\Gamma $ be} a {3RDF} with weight {at most} $k$. Let $%
z=(V_{0},V_{1},V_{2},V_{3},V_{4})$ be one that {label}s large values to $%
w_{j}$'s and $y_{i}$'s. By this choice of $z,$ it is clear that for {any} $j$%
, vertex $w_{j}$ is {label}ed a $4$ and its leaves are {label}ed a $0.$
Observe that vertex $c_{j}$ may also be {label}ed a $0$ under $z.$ Hence the
total weight {label}ed to the vertices of $H_{j}$'s is {at least} $4t.$
Moreover, $z(x_{i})+z(y_{i})\geq 3 $ for every $i,$ and since $k=4t+11q$, we
deduce that $z(x_{i})+z(y_{i})\leq 4.$ We note that by the definition of a {%
3RDF} and our choice of $z$, $z(y_{i})\in \{3,4\}$ and $z(x_{i})=0$ for
every $i.$ We also note that if {$z(y_{i})=3$} for some $i,$ then $x_{i}$
needs that $z(c_{r})=2$ for some $c_{r}\in N(x_{i}).$ {Now let} $\left\vert
\{y_{1,}y_{2},\ldots,y_{3q}\}\cap V_{3}\right\vert =a$ and $\left\vert Z\cap
V_{2}\right\vert =b.$ {Clearly}, $a\leq 3q.$ Also, since every $c_{j}$ has
exactly three {neighbour}s in $X,$ we have $3b\geq a.$ Now, since $%
z(V(G))=4t+2b+3a+4(3q-a)\leq k=4t+11q,$ we deduce that $a-2b\geq q.$
Combining the previous three inequalities we arrive at $b=q$ and $a=3q.$
Consequently, $C^{\prime }=\{C_{j}:z(c_{j})=2\}$ is an exact cover for $C$. $%
\ \Box $ \medskip

\bigskip

Now we can consider the same graph $\Gamma$ built for the transformation in
the proof of Theorem \ref{complexity_bipartite} and add all edges between
the $c_{j}$'s to obtain a chordal graph $\Gamma^{\prime}$. Therefore the
next result is obtained by using the same proof as before on the graph $%
\Gamma^{\prime}.$

\begin{thm}
Problem TRIPLE ROM-DOM is NP-Complete for chordal graphs.
\end{thm}

{\ Finally, we will show that we can solve TRIPLE\ ROM-DOM in linear time
for the class of graphs with bounded clique-width. {Clearly}, this fact
implies the mentioned problem can also be solved in linear time for the
class of trees. }

In what follows, we use several known results related to logic structures.
We refer the reader to the {manuscript}s by Courcelle et al. \cite{CMR} and
by Liedloff et al. \cite{LKLP} for formal details. Namely, we call a $k$%
-expression of $\Gamma,$ on the vertices $\{v_{i}\},$ with labels $%
\{1,2,\ldots ,k\}$ to an expression describing the graph by using the
following operations:

\begin{tabular}{rl}
$\bullet i(x)$ & create a new vertex $x$ with label $i$ \\
$\Gamma_{1}\oplus \Gamma_{2}$ & create a new graph which is the disjoint
union of the graphs $\Gamma_{i}$ \\
$\eta _{ij}(\Gamma)$ & add all edges in $\Gamma$ joining vertices with label
$i$ to vertices with label $j$ \\
$\rho _{i\rightarrow j}(\Gamma)$ & change the label of all vertices with
label $i$ into label $j$%
\end{tabular}

\noindent The \textit{clique-width} of a graph $\Gamma$ is the minimum
integer $k$ which is needed to give a $k$-expression of the graph $\Gamma$.
As an example, we can describe the complete graph $K_3,$ whose set of
vertices is $\{a,b,c\},$ by means of the $2$-expression,
\begin{equation*}
\rho_{2\rightarrow 1}\left( \eta_{12} \left( \rho_{2\rightarrow 1}\left(
\eta_{12} \left( \bullet 1(a) \oplus \bullet 2(b) \right) \right) \oplus
\bullet 2(c) \right) \right)
\end{equation*}
Let us denote by MSOL($\tau_1$) the monadic second order logic with
quantification over subsets of vertices. We also write $\Gamma(\tau_1)$ for
the logic structure $<V(\Gamma),R>,$ where $R$ is a binary relation such
that $R(u,v)$ holds whenever $uv$ is an edge in $\Gamma$.

An optimization problem is said to be a \textit{LinEMSOL($\tau $)
optimization problem} when it is possible to describe it in the following
way (see \cite{LKLP} for more details, since this is a version of the
definition given by \cite{CMR} restricted to finite simple graphs),
\begin{equation*}
\opt\;\left\{ \sum_{1\leq i\leq l}a_{i}|X_{i}|\;:\;<G(\tau
_{1}),X_{1},\ldots ,X_{l}>\;\vDash \theta (X_{1},\ldots ,X_{l})\right\}
\end{equation*}%
where $\theta $ is an MSOL($\tau _{1}$) formula that contains free
set-variables $X_{1},\ldots ,X_{l},$ integers $a_{i}$ and the operator \emph{%
Opt} is either $\min $ or $\max .$

We use the following result on LinEMSOL optimization problems.

\begin{thm}
\emph{(Courcelle et all. \cite{CMR})} Let $k\in \mathbb{N}$ and let $%
\mathcal{C}$ be a class of graphs of clique-width {at most} $k$. Then every
LinEMSOL$(\tau_1)$ optimization problem on $\mathcal{C}$ can be solved in
linear time if a $k$-expression of the graph is part of the input.
\end{thm}

\noindent We extend a result proved by Liedloff et al. (see Th. 31 in \cite%
{LKLP}) regarding the complexity of the {RD-number} decision problem to the
corresponding decision problem for the {TRD-number}.

\begin{thm}
Problem TRIPLE ROM-DOM is a LinEMSOL$(\tau _{1})$ optimization problem.
\end{thm}

\noindent \textbf{Proof. } Let us prove that the TRIPLE ROM-DOM can be
expressed as a LinEMSOL$(\tau _{1})$ optimization problem. Let $%
f=(V_{0},V_{1},V_{2},V_{3},V_{4})$ be a {3RDF} in $\Gamma =(V,E) $ and let
us define the free set-variables $X_{i}$ such that $X_{i}(y)=1$ whenever $%
y\in V_{i}$ and $X_{i}(y)=0,$ {\ for the remaining vertices}. For the sake
of congruence with the logical system notation, we denote by $%
|X_{i}|=\sum_{y\in V}X_{i}(y),$ even when, {clearly}, is $|X_{i}|=|V_{i}|.$

{Note that to solve the decision problem associated with the {TRD-number}
problem is exactly the same to that to reach the optimum for the following
expression.}
\begin{equation*}
\min_{X_i} \; \left\{ \sum_{i=1}^{4}\; i |X_i| \;:\;
<G(\tau_1),X_0,X_1,X_2,X_3,X_4 > \; \vDash \theta(X_0,X_1,X_2,X_3,X_4)
\right\}
\end{equation*}
where $\theta$ is the formula given by
\begin{equation*}
\begin{array}{rl}
\theta(X_0,\ldots,X_4) = & \left( \forall y \left( X_3(y) \lor X_4(y)
\right) \right) \lor \\[.5em]
& \left( X_2(y) \land \exists z \left( R(z,y) \land \left( X_2(z)\lor X_3(z)
\lor X_4(z) \right) \right)\right) \lor \\[.5em]
& \left( X_1(y) \land \exists z \left( R(z,y) \land \left( X_3(z) \lor
X_4(z) \right) \right)\right) \lor \\[.5em]
& \left( X_1(y) \land \exists z,t \left( R(z,y) \land R(t,y) \land X_2(z)
\land X_2(t) \right) \right) \lor \\[.5em]
& \left( X_0(y) \land \exists z \left( R(z,y) \land X_4(z) \right) \right)
\lor \\[.5em]
& \left( X_0(y) \land \exists z,t \left( R(z,y) \land R(t,y) \land X_2(z)
\land X_3(t) \right) \right) \lor \\[.5em]
& \left( X_0(y) \land \exists z,t,s \left( R(z,y) \land R(t,y) \land R(s,y)
\land X_2(z) \land X_2(t) \land X_2(s) \right) \right)%
\end{array}%
\end{equation*}

It is straightforward to verify that $\theta (X_{i})$ is an MSOL$(\tau _{1})$
formula that corresponds to the conditions required for a labeling of the
vertices of the graph to be a {TRD-number} assignment. Namely, the formula
consists of seven clauses and one of which, {at least}, must be true. The
first clause of the formula verifies whether a vertex has a label $3$ or a
label $4$, in which case, no additional condition has to be demanded. In
case of the vertex has a label $2$, the second clause checks if it has a {%
neighbour} with a label $2,3$ or $4,$ and so on. Hence, when the formula $%
\theta (X_{i})$ is satisfied, the requirements of a {3RDF} in $G$ {occurs}. $%
\ \Box $

As a consequence, we may derive the following corollary

\begin{cor}
Problem TRIPLE\ ROM-DOM can be solved in linear time on any graph $\Gamma$
with clique-width bounded by a constant $k$, provided that either there
exists a linear-time algorithm to construct a $k$-expression of $\Gamma$, or
a $k$-expression of $\Gamma$ is part of the input.
\end{cor}

Since any graph with bounded treewidth is also a bounded clique-width graph,
and it is well-known that any tree graph has treewidth equal to 1, then we
can deduce that the TRIPLE ROM-DOM can be solved in linear time for the
class of trees. Besides, there are several classes of graphs $\Gamma$ with
bounded clique-width $cw(\Gamma)$ like, for example, the cographs ($%
cw(\Gamma)\leq 2$) and the distance hereditary graphs ($cw(\Gamma)\leq 3$),
for which it is also possible to solve TRIPLE ROM-DOM in linear time.

\section{Bounds in terms of $p,\Delta $ and $\protect\delta $}

The purpose of this section is to provide various upper and lower bounds on
the {TRD-number} in terms of the order, maximum and minimum degrees of a
graph. Since the function that {label}s $2$ to each vertex in {an ntc-graph}
$\Gamma$ is a 3RDF, the following observation is immediate.

\begin{obs}
\label{cota2n}For any {ntc-graph} the inequality $\gamma _{\lbrack
3R]}\left( {\Gamma}\right) \leq 2p$ holds.
\end{obs}

Through the following results, we prove upper bounds on the {TRD-number} of {%
an ntc-graph} improving, in some cases, the one given in Observation \ref%
{cota2n}.

\begin{prop}
\label{cota1}Let $\Gamma$ be {an ntc-graph} with $p\geq 2$ vertices and
maximum degree $\Delta \geq 1$. Then $\gamma _{\lbrack 3R]}\left(
\Gamma\right) \leq 3p-3\Delta +1.$
\end{prop}

\noindent \textbf{Proof.} Consider a vertex $v\in V(\Gamma)$ of maximum
degree $\Delta $ and define the function $h:V\rightarrow \{0,1,2,3,4\}$ as
follows: $h(v)=4,h(w)=0$ for all $w\in N(v)$ and $h(w)=3$ {\ for the
remaining vertices}. It is straightforward to see that $h$ is a 3RDF, and
hence
\begin{equation*}
\gamma _{\lbrack 3R]}\left( \Gamma\right) \leq h(V)=4+3\left( p-1-\Delta
\right) =3p-3\Delta +1.
\end{equation*}%
$\ \Box $

\noindent Note that for any {ntc-graph} satisfying that $\Delta (\Gamma)>%
\frac{p+1}{3}$, the upper bound given in Proposition \ref{cota1} is better
than the more general upper bound pointed out in Observation \ref{cota2n}.

\begin{cor}
Let $\Gamma$ be {an ntc-graph} with $p\geq 2$ vertices, {$girth \ge 4$,}
maximum degree $\Delta \le p-2$ and minimum degree $\delta \geq 2$. Then $%
\gamma _{\lbrack 3R]}\left( \Gamma\right) \leq 3p-3\Delta .$
\end{cor}

\noindent \textbf{Proof.} {I}t is straightforward to check that the function
defined in the proof of Proposition \ref{cota1} can be modified such that $%
h(v)=3$ and, since $\delta \geq 2$ {and $g(\Gamma )\geq 4$}, it is still a
3RDF. $\ \Box $

\noindent Next, we prove that the previous upper bound can be improved if
the graph meets certain requirements.

\begin{prop}
\label{cota1bis}Let $\Gamma$ be {an ntc-graph} with $p\geq 2$ vertices,
minimum degree $\delta (\Gamma)\geq 2$ and girth {at least} $5$. Then $%
\gamma _{\lbrack 3R]}\left( \Gamma\right) \leq 2p-2\Delta (\Gamma)+1.$
\end{prop}

\noindent \textbf{Proof.} Let us denote by $v$ a vertex with maximum degree $%
d(v)=\Delta (\Gamma)$. Consider the function $l:V(\Gamma)\rightarrow
\{0,1,2,3,4\}$ defined as follows: $l(v)=3$, $l(w)=0$ for all $w\in N(v)$
and $l(w)=2$ {\ for the remaining vertices}. Since $\delta (\Gamma)\geq 2$
and {$g(\Gamma)\geq 5$}, each vertex in {$V(\Gamma)-\{v\}$} must have {at
least} a {neighbour} in the set $V(\Gamma)-N[v].$ Therefore, any vertex {%
label}ed a $2$ is {connecting to}, {at least}, another vertex with {the same
value. it follows} that $l$ is a 3RDF in $\Gamma$ and, in consequence,
\begin{equation*}
\gamma _{\lbrack 3R]}\left( \Gamma\right) \leq 3+2(p-\Delta -1)=2p-2\Delta
+1.
\end{equation*}%
$\ \Box $

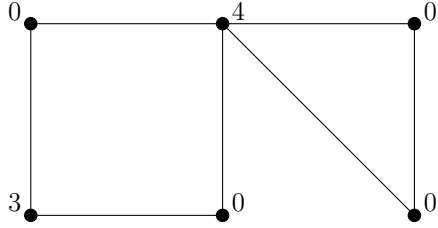
\begin{figure}[h]
\centering
\begin{tikzpicture}[scale=.85, transform shape]


\node [draw, shape=circle,fill=black,scale=0.5] (u0) at  (3,2) {};
\node [draw, shape=circle,fill=black,scale=0.5] (u1) at  (3,5) {};
\node [draw, shape=circle,fill=black,scale=0.5] (u2) at  (6,5) {};
\node [draw, shape=circle,fill=black,scale=0.5] (u3) at  (6,2) {};
\node [draw, shape=circle,fill=black,scale=0.5] (u4) at  (9,5) {};
\node [draw, shape=circle,fill=black,scale=0.5] (u5) at  (9,2) {};
\draw (u0)--(u1);
\draw (u1)--(u2);
\draw (u2)--(u3);
\draw (u3)--(u0);
\draw (u2)--(u4);
\draw (u4)--(u5);
\draw (u2)--(u5);
\node[draw=none, fill=none] at (2.75,2.2){ $3$};
\node[draw=none, fill=none] at (6.25,2.2){ $0$};
\node[draw=none, fill=none] at (2.75,5.2){ $0$};
\node[draw=none, fill=none] at (6.25,5.2){ $4$};
\node[draw=none, fill=none] at (9.25,2.2){ $0$};
\node[draw=none, fill=none] at (9.25,5.2){ $0$};

\end{tikzpicture}
\caption{{\protect\small The condition {$g(\Gamma)\ge 5$} is necessary in
Proposition \protect\ref{cota2} .}}
\label{c4c3}
\end{figure}

It is worth noting that the condition required in Proposition \ref{cota1bis}
about the minimum degree is essential. As an example, we can consider the
path $P_{4},$ where $\gamma _{\lbrack 3R]}\left( P_{4}\right) =7$ which
implies that the upper bound given by Proposition \ref{cota1} is {tight}.
However, since $2p-2\Delta +1=5,$ then for this particular graph the bound
proved in Proposition \ref{cota1bis} does not apply because $\delta
(\Gamma)=1$. Similarly, the condition about the girth of the graph is also
essential, as we may observe by the cycle $C_{4}$ for which $\gamma
_{\lbrack 3R]}\left( C_{4}\right) =6>2p-2\Delta +1=5,$ or the graph depicted
in Figure \ref{c4c3}. The graph has girth $3$ and, in this case, $\gamma
_{\lbrack 3R]}\left( \Gamma\right) =7$ which shows that the bound of
Proposition \ref{cota1} is {tight}, while by applying Proposition \ref%
{cota1bis}, it would be $\gamma _{\lbrack 3R]}\left( \Gamma\right) \leq 5$.

\begin{cor}
\label{cota2}Let $\Gamma$ be a connected $r$-regular graph, with girth {$%
g(\Gamma)\geq 7 $} and order $p$. Then $\gamma _{\lbrack 3R]}\left(
\Gamma\right) \leq 2p-2r^2+3r-2.$
\end{cor}

\noindent \textbf{Proof.} Let $v$ be any vertex of $\Gamma$ and let us
consider $v $ as the root of {a spanning} tree in $\Gamma$. So that, $%
N_0=\{v\}$ is the root, $N_1=N(v)$ is the first level and $|N_1(v)|=r$, $%
N_2=N(N_1)-N_0$ is the second level and $|N_2(v)|=r(r-1)$ and so on. Define
the function $h:V\rightarrow \{0,1,2,3,4\}$ as follows: $h(w)=3$ for all $%
w\in N_1$, $h(w)=0$ for all $w\in N_0\cup N_2$ and $h(w)=2$ {\ for the
remaining vertices}. {Clearly}, $h$ is a 3RDF and the weight of $h$ is:
\begin{equation*}
h(V)=3r+2(p-1-r-r(r-1))
\end{equation*}%
which is well-defined since {$g(\Gamma)\geq 7$}, thus $p\geq r^{3}-r^2+r+1,$
and every vertex in $N_3\cap V_2$ must have a neighbour in $V_2.$ Therefore,
\begin{equation*}
\gamma _{\lbrack 3R]}\left( \Gamma\right) \leq
h(V)=3r+2(p-r^2-1)=2p-2r^{2}+3r-2.
\end{equation*}%
$\ \Box $

Observe that $2r^{2}-3r+2\ge 2r-1$ for all $r\ge 1$ and therefore the upper
bound given by Corollary \ref{cota2} is better than the one given in
Proposition \ref{cota1bis} for $r$-regular graphs with girth {at least} $7$.
Besides, the upper bound given in Corollary \ref{cota2} is {tight} for $C_7,$
among other graphs.

Using a probabilistic method, we provide the next upper bound on $\gamma
_{\lbrack 3R]}\left( G\right) .$

\begin{prop}
\label{prop-probab}Let $\Gamma$ be a graph of order $p$, minimum degree $%
\delta $ and maximum degree $\Delta $. Then
\begin{equation*}
\gamma _{\lbrack 3R]}\left( \Gamma\right) \leq \left\lfloor\frac{4p}{\delta
+1}\left( \ln \left( \frac{3(\delta +1)}{4}\right) +1\right)\right\rfloor .
\end{equation*}
\end{prop}

\noindent \textbf{Proof.} Let us denote by $p^{\prime}\in (0,1)$ the
probability that any vertex $v$ belongs to certain set of vertices $%
A\subseteq V(\Gamma).$ Let us denote by $B=N[A]^{c}=A^{c}\cap N(A)^{c}.$
This means that
\begin{equation*}
\begin{array}{rcl}
P[v\in B] & = & P[\mbox{neither v nor its neighbours belong to A}] \\%
[0.75em]
& = & (1-p^{\prime})(1-p^{\prime})^{d(v)} \\[0.75em]
& = & (1-p^{\prime})^{d(v)+1}\leq (1-p^{\prime})^{\delta (\Gamma)+1},%
\end{array}%
\end{equation*}%
because $0<p^{\prime}<1.$ {Clearly}, the expected size of the set $A$ is
equal to $E[|A|]=pp^{\prime}$ and the corresponding expected value of $|B|$
is $E[|B|]\leq p(1-p^{\prime})^{\delta (\Gamma)+1}.$ Now, let us define the
following function on $V(\Gamma),$
\begin{equation*}
f(v)=\left\{
\begin{array}{ll}
4 & \mbox{if $v\in A$,} \\[0.5em]
0 & \mbox{if $v\in N(A)$,} \\[0.5em]
3 & \mbox{if $v\in B=V(\Gamma)-N[A].$}%
\end{array}%
\right.
\end{equation*}%
The expected value of $f(V)$ would be
\begin{equation*}
E[f(V)]=4E[|A|]+3E[|B|]\leq 4pp^{\prime}+{3p}(1-p^{\prime})^{\delta
(\Gamma)+1}
\end{equation*}%
Taking into account that $(1-p^{\prime})<\mathrm{e}^{-p^{\prime}}$ when $%
0<p^{\prime}<1,$ it is deduced that
\begin{equation*}
E[f(V)]\leq 4pp^{\prime}+3p\mathrm{e}^{-p^{\prime}(\delta (\Gamma)+1)}
\end{equation*}%
The value $p^{\prime}$ for which the minimum value of the latter expression
is reached satisfies $4p-3p(\delta (\Gamma)+1)\mathrm{e}^{-p^{\prime}(\delta
(\Gamma)+1)}=0$, {that} implies that
\begin{equation*}
\mathrm{e}^{-p^{\prime}(\delta (\Gamma)+1)}=\frac{4}{3(\delta (\Gamma)+1)}
\end{equation*}%
This leads us to the desired value of the probability $p^{\prime}=\frac{1}{%
\delta (\Gamma)+1}\ln \left( \frac{3(\delta (\Gamma)+1)}{4}\right) ,$ and
finally
\begin{equation*}
\gamma _{\lbrack 3R]}\left( \Gamma\right) \leq 4p\frac{1}{\delta (\Gamma)+1}%
\ln \left( \frac{3(\delta (\Gamma)+1)}{4}\right) +3p\frac{4}{3(\delta
(\Gamma)+1)}
\end{equation*}%
which concludes the proof. $\ \Box $

It is worth noting that the upper bound provided by Proposition \ref%
{prop-probab} may be smaller or larger than the one given by Proposition \ref%
{cota1}. To see consider the cycle $C_{5},$ where by Proposition \ref{cota1}%
, we obtain that $\gamma _{\lbrack 3R]}(C_{5})\leq 10$ which is better than
the bound derived from the probabilistic procedure, $\gamma _{\lbrack
3R]}(C_{5})\leq 12.$ On the contrary, for the cycle $C_{9},$ Proposition \ref%
{prop-probab} provides $\gamma _{\lbrack 3R]}(C_{9})\leq 21$ while by
Proposition \ref{cota1}, we have $\gamma _{\lbrack 3R]}(C_{9})\leq 22.$

Our next result shows that any {ntc-graph} $\Gamma$ with order $p\geq 2$, $%
\gamma _{\lbrack 3R]}(\Gamma)\leq \frac{7p}{4}$. {Since removing} an edge
cannot decrease the {TRD-number}, it suffices to prove the result for trees.
We then characterize the {ntc-graph}s attaining the upper bound.

Let $\mathcal{F}$ be the family of all trees that can be built from {$k\geq 1
$ paths $P_{4}^{i}:=v_{1}^{i}v_{2}^{i}v_{3}^{i}v_{4}^{i}\;(1\leq i\leq k)$
by adding $k-1$ edges incident with the $v_{2}^{i}$'s so that they induce a
connected subgraph.} We write \textquotedblleft Ind-Hyp\textquotedblright\
instead of \textquotedblleft induction hypothesis\textquotedblright .

\begin{prop}
\label{tree1} \emph{If $T\in \mathcal{F}$, then $\gamma_{[3R]}(T)\ge \frac{%
7p(T)}{4}$.}
\end{prop}

\noindent \textbf{Proof.} Suppose $T$ is a tree in $\mathcal{F}$ and let $T$
be built from $k$ copies of $P_{4}$. Note that if $H$ is an induced subgraph
of $T$ {for which} $H\cong P_{4}$ and the leaves of $H$ are leaves of $T$,
then every 3RDF of $T$ {label}s a weight of {at least} 7 to $H$. Since $T$
has $p(T)/4$ disjoint copies of the induced subgraph $P_{4}$, where the
leaves of the $P_{4}$ are leaves of $T$, we have $\gamma _{\lbrack
3R]}(T)\geq \frac{7p(T)}{4}$. $\hfill \Box $

\begin{thm}
\label{main} \emph{If $T$ is a tree with order $p\ge 3$, then $%
\gamma_{[3R]}(T)\le \frac{7p}{4}$.}
\end{thm}

\noindent \textbf{Proof.} The proof is by induction on $p$. Since $p\geq 3$,
we have $\mathrm{diam}(T)\geq 2$. If $\mathrm{diam}(T)=2$, then $T$ is a
star and we have $\gamma _{\lbrack 3R]}(T)=4<\frac{7p}{4}$. If the diameter
of $T$ is 3, then $T$ is a double star {$DS_{r,s}$ for $r\geq s\geq 1$. If $%
s=1$, then $\gamma _{\lbrack 3R]}(T)=7\leq \frac{7(3+r)}{4}$ with equality
if and only if {$r=1,$ that is $T=P_{4}$. If $s\geq 2$,} then $%
p\geq 6$} and $\gamma _{\lbrack 3R]}(T)=8<\frac{7p}{4}$. Hence, we may {%
suppose} that $\mathrm{diam}(T)\geq 4$ and this implies that {$p\geq 5$.} {%
Assume} that the statement is true for {each} tree $T^{\prime }$ of order $%
p^{\prime }$ with $3\leq p^{\prime }<p$. Let $T$ be a tree of order $p$ and
let $v_{1}v_{2}\ldots v_{k}$ be a diametral path in $T$ such that $d(v_{2})$
is as large as possible. {\color{green} Assume, without loss of generality, that
the tree $T$ is rooted at the vertex $v_k$.} {For sake of
notation, we will often write }$T-T_{v_{j}}$ {to denote the tree obtained
from }$T$ {by removing }$v_{j}$ and all its descendants. If $d(v_{2})\geq 3$%
, then any $\gamma _{\lbrack 3R]}(T-T_{v_{2}})$-function $h$ can be extended
to a 3RDF of $T$ {by keeping the assignments given to vertices of }$%
T-T_{v_{2}}$ {under }$h$ {to which we also }label $v_{2}$ by 4 and the
leaves {in }$L_{v_{2}}$ by 0. {It }follows from the {Ind-Hyp} that
\begin{equation*}
\gamma _{\lbrack 3R]}(T)\leq \gamma _{\lbrack 3R]}(T-T_{v_{2}})+4\leq \frac{%
7(p-3)}{4}+4<\frac{7p}{4}.
\end{equation*}%
Hence assume that $d(v_{2})=2$. By the choice of diametral path, we may {%
suppose} that any child of $v_{3}$ with depth one is of degree 2. If $%
d(v_{3})=2$ and $T-T_{v_{3}}=P_{2}$, then $T=P_{5}$ and {clearly} $\gamma
_{\lbrack 3R]}(T)<\frac{7p}{4}$ and if $d(v_{3})=2$ and $p(T-T_{v_{3}})\geq
3 $, then any $\gamma _{\lbrack 3R]}(T-T_{v_{3}})$-function can be extended
to a 3RDF of $T$ by {labeling} a 4 to $v_{2}$ and a 0 to $v_{1},v_{3}$ and
as above we have $\gamma _{\lbrack 3R]}(T)<\frac{7p}{4}.$ Henceforth we {%
assume} that $d(v_{3})\geq 3$. We consider the following cases.\newline
\newline
\smallskip \noindent \textbf{Case 1.} $d(v_{3})\geq 4$.\newline
Let $t_{1}$ be the number of children of $v_{3}$ with depth one and $%
t_{2}=|L_{v_{3}}|$, $t_{1}+t_{2}\geq 3$ and $t_{1}\geq 1$. Set $t=0$ if $%
t_{2}=0$ and $t=1$ if $t_{2}\geq 1$. Let $T^{\prime }=T-T_{v_{3}}$. If {%
$p(T^{\prime })=2$}, then $T$ is a spider and {labeling} $3+t$
to $v_{3}$ and a 3 to each leaf at distance two from $v_{3}$ provides a 3RDF
of $T$ of weight $3t_{1}+3+t$ and this implies that $\gamma _{\lbrack
3R]}(T)\leq 3t_{1}+3+t<\frac{7(2t_{1}+1+t_{2})}{4}=\frac{7p}{4}.$ Let $%
p(T^{\prime })\geq 3$. Then any $\gamma _{\lbrack 3R]}(T-T_{v_{3}})$%
-function can be extended to a 3RDF of $T$ by {labeling} $3+t$ to $v_{3}$, a
3 to all leaves of $T_{v_{3}}$ at distance two from $v_{3}$ and 0 to other
vertices and by the {Ind-Hyp} we have
\begin{equation*}
\gamma _{\lbrack 3R]}(T)\leq \gamma _{\lbrack
3R]}(T-T_{v_{3}})+3t_{1}+3+t\leq \frac{7(p-(2t_{1}+1+t_{2}))}{4}+3t_{1}+3+t<%
\frac{7p}{4}.
\end{equation*}

\smallskip \noindent \textbf{Case 2.} $d(v_{3})=3$ and $v_{3}$ is not a {\
stem} vertex.\newline
Then $T_{v_{3}}$ is a healthy spider with two feet. The result is immediate
if $T-T_{v_{3}}=P_{2}$. So assume that $p(T-T_{v_{2}})\geq 3$. Then any $%
\gamma _{\lbrack 3R]}(T-T_{v_{3}})$-function can be extended to a {3RDF} of $%
T$ by {labeling} $4$ to children of $v_{3}$ and a 0 to $v_{3}$ and all
leaves of $T_{v_{3}}$ and by the {Ind-Hyp} we have
\begin{equation*}
\gamma _{\lbrack 3R]}(T)\leq \gamma _{\lbrack 3R]}(T-T_{v_{3}})+8\leq \frac{%
7(p-5)}{4}+8<\frac{7p}{4}.
\end{equation*}

\smallskip \noindent \textbf{Case 3.} $d(v_{3})=3$ and $v_{3}$ is a {\ stem}
vertex.\newline
Let $w$ be the leaf {connecting to} $v_{3}$. Then $T_{v_{3}}\cong P_{4}.$ If
$p=6,$ then {clearly} $\gamma _{\lbrack 3R]}(T)\leq 10<\frac{7p}{4}.$ Hence
let $p\geq 7.$ Then any $\gamma _{\lbrack 3R]}(T-T_{v_{3}})$-function can be
extended to a {3RDF} of $T$ by {labeling} $4$ to $v_{3}$, a 3 to $v_{1}$ and
a 0 to $v_{2},w$ and by the {Ind-Hyp} we obtain
\begin{equation}
\gamma _{\lbrack 3R]}(T)\leq \gamma _{\lbrack 3R]}(T-T_{v_{3}})+7\leq \frac{%
7(p-4)}{4}+7=\frac{7p}{4}.  \label{pop}
\end{equation}%
$\hfill \Box $

\begin{thm}
\emph{Let $T$ be a tree of order $p\ge 3$. Then $\gamma_{[3R]}(T)=\frac{7p}{4%
}$ if and only if $T\in \mathcal{F}$.}
\end{thm}

\noindent \textbf{Proof.} The sufficiency follows from Proposition \ref%
{tree1} and Theorem \ref{main}. To show that every tree $T$ with $\gamma
_{\lbrack 3R]}(T)=\frac{7p}{4}$ is in $\mathcal{F}$, we {demand to the}
proof of Theorem \ref{main}. Assume that $\gamma _{\lbrack 3R]}(T)=\frac{7p}{%
4}$. The proof is by the induction on $p$. Since $P_{4}\in \mathcal{F}$, and
$p$ is a multiple of $4$, we may assume that $p\geq 8$. Following the proof
of Theorem \ref{main}, there is only one case, namely, Case 3, where it is
possible to achieve equality. Using the terminology from this proof, we have
$\gamma _{\lbrack 3R]}(T-T_{v_{3}})=\frac{7p(T-T_{v_{3}})}{4}$. It follows
from the {Ind-Hyp} that $T-T_{v_{3}}\in \mathcal{F}$. Thus $T-T_{v_{3}}$ is
built from {$l=p/4-1$ paths $P_{4}^{i}:=v_{1}^{i}v_{2}^{i}v_{3}^{i}v_{4}^{i}%
\;(1\leq i\leq \ell )$ by adding $\ell -1$ edges incident to the vertices $%
v_{2}^{1},v_{2}^{2},\ldots ,v_{2}^{\ell }$} {so that they induce a subtree}.
{Hence vertex }$v_{4}$ {belongs to some path }{$P_{4}^{i},$} {that we
denote, without loss of generality, by }{$H=P_{4}^{\ell },$} constituting
the tree $T-T_{v_{3}}.$ {Note that any $\gamma _{\lbrack 3R]}(T-T_{v_{3}})$%
-function labels a weight of 7 to each path $P_{4}^{i}\;(1\leq i\leq \ell )$.%
} {We claim that $v_4$ is a stem of $H=P_{4}^{\ell }$. Suppose, to the contrary, that $v_4$ is a leaf in $H=P_{4}^{\ell }$.
If $v_{4}=v_{4}^{\ell }$, then the
function $z:V(T)\rightarrow \{0,1,2,3,4\}$ defined by {$z(v_{1})=3,z(v_{3})=4
$,} $z(v_{2}^{i})=4$ for $1\leq i\leq \ell $, $z(v_{4}^{i})=3$ for each $%
1\leq i\leq \ell -1$ and $z(x)=0$ for the remaining vertices, is a 3RDF of $%
T$ of weight {$\frac{7p}{4}-3$} which leads to a contradiction. If $v_{4}=v_{1}^{\ell }$, then the
function $z:V(T)\rightarrow \{0,1,2,3,4\}$ defined by $z(v_{1})=3,z(v_{3})=4, z(v_{3}^{\ell})=4$,
$z(v_{2}^{i})=4$ for $1\leq i\leq \ell-1 $, $z(v_{4}^{i})=3$ for each $%
1\leq i\leq \ell -1$ and $z(x)=0$ for the remaining vertices, is a 3RDF of $%
T$ of weight {$\frac{7p}{4}-3$} which is a contradiction again.}
Thus $v_{4}$ is a stem vertex of {$H=P_{4}^{\ell }$. If $\ell =1$, then
clearly $T\in \mathcal{F}$. Assume that $\ell \geq 2$. If $v_{4}=v_{3}^{\ell
}$, then the function $z:V(T)\rightarrow \{0,1,2,3,4\}$ defined by $%
z(v_{1})=3,z(v_{3})=4$, $z(v_{1}^{\ell })=z(v_{4}^{\ell })=3$, $%
z(v_{2}^{i})=4$ for $1\leq i\leq \ell -1$, $z(v_{4}^{i})=3$ for $1\leq i\leq
\ell -1$ and $z(x)=0$ for the remaining vertices, is a 3RDF of $T$ of weight
{$\frac{7p}{4}-3$} which leads to a contradiction. Thus $%
v_{4}=v_{2}^{\ell }$} and this implies that $T\in \mathcal{F}$. $\hfill \Box
$

The next result is an immediate consequence of Theorem \ref{main}.

\begin{cor}
\emph{If $\Gamma$ is {an ntc-graph} with order $p\ge 3$, then $%
\gamma_{[3R]}(\Gamma)\le \frac{7p}{4}$.}
\end{cor}

{\ Assume $\mathcal{H}$ is the family of {ntc-graph}s $\Gamma $ of order $p$
that {can be built from $\ell =p/4$ paths $%
P_{4}:=v_{1}^{i}v_{2}^{i}v_{3}^{i}v_{4}^{i}\;(1\leq i\leq \ell )$ by adding
some edges between the vertices $v_{2}^{1},v_{2}^{2},\ldots ,v_{2}^{\ell }$
so that they induce a connected subgraph.} %
}

\begin{thm}
\emph{The family $\mathcal{H}$ is precisely the family of {ntc-graph}s $%
\Gamma$ {such that} $\gamma_{[3R]}(\Gamma)= \frac{7p}{4}$.}
\end{thm}

\noindent \textbf{Proof.} If $\Gamma \in \mathcal{H}$, then as {we did for
the family }$\mathcal{F}$,{\ }we can see that $\gamma _{\lbrack 3R]}(\Gamma
)=\frac{7p}{4}$. Now assume that $\gamma _{\lbrack 3R]}(\Gamma )=\frac{7p}{4}
$. {Since removing} edges cannot decrease $\gamma _{\lbrack 3R]}(\Gamma )$,
we deduce that {any} spanning tree of $\Gamma $ has {TRD-number} $\frac{7p}{4%
}$. Hence, {every} spanning tree of $\Gamma $ is in $\mathcal{F}$. {In fact
since any spanning tree $T$ {of $\Gamma $ can be built from $\ell =p/4$
paths $P_{4}^{i}:=v_{1}^{i}v_{2}^{i}v_{3}^{i}v_{4}^{i}\;(1\leq i\leq \ell )$
by adding $\ell -1$ edges incident with vertices $v_{2}^{1},v_{2}^{2},\ldots
,v_{2}^{\ell }$ so that they induce a subtree, any edge that does not belong
to }$T$ has its endvertices in $\{${$v_{2}^{1},v_{2}^{2},\ldots ,v_{2}^{\ell
}\}$, and therefore}, $\Gamma \in \mathcal{H}$.}


\section{Relationships with some Roman domination parameters}

In 2004, Cockayne et al. showed that for any graph, $\gamma (\Gamma)\leq
\gamma _{R}(\Gamma)\leq 2\gamma (\Gamma),$ (see \cite{CDHH}, Prop. 1)$.$ In
2016, Beeler et al. proved that $2\gamma (\Gamma)\leq \gamma _{dR}(\Gamma),$
(see \cite{BHH}, Prop. 8). Moreover, in the same {manuscript}, it is proved
that $\gamma _{dR}(\Gamma)\leq 2\gamma _{R}(\Gamma).$ Next, we give
relations involving $\gamma _{\lbrack 3R]}\left( \Gamma\right) $ with $%
\gamma (\Gamma),\gamma _{R}$ and $\gamma _{dR}\left( \Gamma\right) $ that
complete known chain of inequalities.

\begin{prop}
\label{nounos}For any {ntc-graph} $\Gamma ,$ there exists a $\gamma
_{\lbrack 3R]}\left( \Gamma \right) $-function that does not {assign a }$1$
to any vertex in $\Gamma $.
\end{prop}

\noindent \textbf{Proof.} Among all $\gamma _{\lbrack 3R]}\left(
\Gamma\right) $-functions let $h=\left( V_{0},V_{1},V_{2},V_{3},V_{4}\right)
$ be one such that $\left\vert V_{1}\right\vert $ is as small as possible.
If $V_{1}=\emptyset ,$ then $h$ is the desired function. Hence we assume
that $V_{1}\neq \emptyset ,$ and let $v$ be a vertex of $\Gamma$ such that $%
h(v)=1.$ Since $h$ is a 3RDF then either there is a vertex $w\in
N_{\Gamma}(v)$ with $h(w)\geq 3$ or there are two vertices $w_{1},w_{2}\in
N(v)\cap V_{2}.$ In the first case, taking into account that $h$ has minimum
weight, it is clear that $h(w)=3$. In the former case, the new function $l$
defined as follows, $l(u)=h(u)$ for every $u\in V(\Gamma)-\{v,w\}$ and $%
l(v)=0,l(w)=4$ is a $\gamma _{\lbrack 3R]}\left( \Gamma\right) $-function.
In the latter case, the function $z$ defined by $z(u)=h(u)$ for every $u\in
V(\Gamma)-\{v,w_{1}\},$ $l(v)=0,$ $l(w_{1})=3$ is a $\gamma _{\lbrack
3R]}\left( \Gamma\right) $-function. In either case, $l$ or $z$ {label}s $1$%
s to vertices less than $h,$ a contradiction. $\Box $

\begin{prop}
\label{cotageneral}Let $\Gamma$ be any {ntc-graph}. Then
\begin{equation*}
\gamma (\Gamma)\leq \gamma _{R}(\Gamma)\leq 2\gamma (\Gamma)\leq \gamma
_{dR}(\Gamma)<\gamma _{\lbrack 3R]}\left( \Gamma\right) \leq \min \{\frac{3}{%
2}\gamma _{dR}(\Gamma),4\gamma (\Gamma)\}
\end{equation*}
\end{prop}

\noindent \textbf{Proof.} As mentioned before, it is sufficient to prove
that $\gamma _{dR}(\Gamma )<\gamma _{\lbrack 3R]}\left( G\right) \leq \min \{%
\frac{3}{2}\gamma _{dR}(\Gamma ),4\gamma (\Gamma )\}.$ By Proposition \ref%
{nounos}, there is a $\gamma _{\lbrack 3R]}\left( \Gamma \right) $-function
with no vertex {label}ed a $1.$ So, let $h=(V_{0},\emptyset
,V_{2},V_{3},V_{4})$ be such a {3RDF} of minimum weight, and let us define
the function $l$ as follows: $l(v)=3$ for all $v\in V_{3}\cup V_{4}$ and $%
l(v)=h(v)$ {\ for the remaining vertices}. {Clearly}, $l$ is a {DRDF} in $%
\Gamma $ and
\begin{equation*}
\gamma _{dR}(\Gamma )\leq l(V)=h(V)-|V_{4}|\leq h(V)=\gamma _{\lbrack
3R]}\left( \Gamma \right) .
\end{equation*}%
From the last expression we may derive that $\gamma _{dR}(\Gamma )=\gamma
_{\lbrack 3R]}\left( \Gamma \right) $ implies that $V_{4}=\emptyset ,$ and
therefore $V_{3}\neq \emptyset .$ So, there exists a $\gamma _{\lbrack
3R]}\left( \Gamma \right) $-function $h=(V_{0},\emptyset
,V_{2},V_{3},\emptyset )$ such that $l=(V_{0},\emptyset ,V_{2},V_{3})$ is
also a $\gamma _{dR}$-function in $\Gamma $. Let $v$ be any vertex of $V_{3}$
and let $z$ be the function defined by $z(v)=2$ and $z(w)=l(w)$ for all $%
w\neq v.$ Taking into account that for any vertex $u\in N(v)\cap V_{0}$ it
must be $l(N[u])\geq 5$ we have that $z(N[u])\geq 4.$ Hence $z$ is a {DRDF}
in $\Gamma $ of weight $z(V)=l(V)-1$, which is a contradiction. {Thus}, $%
\gamma _{dR}(\Gamma )<\gamma _{\lbrack 3R]}\left( \Gamma \right) .$

{To} prove the remaining inequality, we first consider a $\gamma (\Gamma)$%
-set $D$, where each of $D$ is {label}ed a 4 and each vertex not in $D$ is {%
label}ed a $0.$ {Clearly} we obtain a 3RDF of $\Gamma$ and so $\gamma
_{\lbrack 3R]}(\Gamma)\leq 4\gamma (\Gamma)$.

Now, to see that $\gamma _{\lbrack 3R]}\left( \Gamma\right) \leq \frac{3}{2}%
\gamma _{dR}(\Gamma),$ let us consider a $\gamma _{dR}(\Gamma)$-function
with no 1 labels, $h=(V_{0},V_{2},V_{3}).$ Note that such a function $h$
exists (see \cite{BHH}). Let $l$ be a function defined as: $l(v)=4$ for all $%
v\in V_{3}$, $l(v)=3$ for all $v\in V_{2}$ and $l(v)=0$ {\ for the remaining
vertices}. Then $l$ is a 3RDF on $V(\Gamma),$ {and so}
\begin{equation*}
\gamma _{\lbrack 3R]}\left( \Gamma\right) \leq l(V)=3|V_{2}|+4|V_{3}|\leq
3|V_{2}|+\frac{9}{2}|V_{3}|=\frac{3}{2}\left( 2|V_{2}|+3|V_{3}|\right) =%
\frac{3}{2}\gamma _{dR}(\Gamma).
\end{equation*}%
$\ \Box $

Our next result is a lower bound for the {TRD-number} of a graph in terms of
the order, the maximum degree and the {domination number} of a graph $\Gamma
$.

\begin{prop}
\label{lower}Let $\Gamma$ be {an ntc-graph} with order {$p$}, maximum degree
$\Delta $ and domination number $\gamma =\gamma (\Gamma)$. Then
\begin{equation*}
\gamma _{\lbrack 3R]}\left( \Gamma\right) \geq \left\lceil \frac{2p+{(\Delta
-1)}\gamma (\Gamma)}{\Delta }\right\rceil.
\end{equation*}
{This bound is {tight} for any graph with a universal vertex.}
\end{prop}

\noindent \textbf{Proof.} Let $h=(V_{0},\emptyset ,V_{2},V_{3},V_{4})$ be a $%
\gamma _{\lbrack 3R]}\left( \Gamma\right) $-function in $\Gamma,$ where $%
\gamma _{\lbrack 3R]}\left( \Gamma\right) =2|V_{2}|+3|V_{3}|+4|V_{4}|.$ Let $%
V_{0}=V_{0}^{2}\cup V_{0}^{3}\cup V_{0}^{4}$ be a partition of the set $%
V_{0} $ such that $V_{0}^{4}=V_{0}\cap N(V_{4}),V_{0}^{3}=\left( V_{0}\cap
N(V_{3})\right) -V_{0}^{4},$ and $V_{0}^{2}=V_{0}-\left( V_{0}^{4}\cup
V_{0}^{3}\right) \subseteq N(V_{2}).$ Since the maximum degree is $\Delta ,$
any vertex of $V_{4}$ can be {connecting to} {at most} $\Delta $ vertices in
$V_{0} $. That is to say, $|V_{0}^{4}|\leq \Delta |V_{4}|$. On the other
hand, each vertex in $V_{0}^{3}$ must have {at least} two neighbours is $%
V_{2}\cup V_{3},$ so $2|V_{0}^{3}|\leq |E[V_{2}\cup V_{3},V_{0}^{3}]|\leq
|E[V_{2},V_{0}^{3}]|+ \Delta |V_{3}|.$ Finally, any vertex in $V_{0}^{2}$ is
{connecting to} {at least} three vertices in $V_{2}$, and thus $%
3|V_{0}^{2}|\leq |E[V_{2},V_{0}^{2}]|$. Summing up the latter bounds we have
that
\begin{equation*}
\begin{array}{rcl}
|V_{0}| & = & |V_{0}^{2}|+|V_{0}^{3}|+|V_{0}^{4}| \\[1 em]
& \leq & \frac{|E[V_{2},V_{0}^{2}]| }{3}+\frac{|E[V_{2},V_{0}^{3}]|}{2} +
\frac{\Delta}{2}|V_{3}| +\Delta |V_{4}| \\[1em]
& \leq & \frac{|E[V_{2},V_{0}^{2}]| + |E[V_{2},V_{0}^{3}]|}{2} + \frac{\Delta%
}{2}|V_{3}| +\Delta |V_{4}| \\[1em]
& \leq & {\frac{\Delta-1 }{2}}|V_{2}|+ \frac{\Delta}{2}|V_{3}| +\Delta
|V_{4}|%
\end{array}%
\end{equation*}%
or, equivalently,
\begin{equation*}
\displaystyle\frac{2}{\Delta }|V_{0}| \leq {\frac{\Delta-1}{\Delta}}%
|V_{2}|+|V_{3}|+2|V_{4}|
\end{equation*}

Now, taking into account that $V_{2}\cup V_{3}\cup V_{4}$ is a dominating
set in $V(\Gamma)$, it is possible to estimate the value of the parameter $%
\gamma_{\lbrack 3R]}\left( \Gamma\right) $ as follows,
\begin{equation*}
\begin{array}{rcl}
\gamma_{\lbrack 3R]}\left( \Gamma\right) & = & 2|V_{2}|+3|V_{3}|+4|V_{4}| \\%
[1em]
& \geq & |V_{2}|+|V_{3}|+|V_{4}|+{\frac{\Delta-1}{\Delta}}%
|V_{2}|+|V_{3}|+2|V_{4}|+{\frac{1}{\Delta}|V_{2}|+|V_{3}|+|V_{4}|} \\[1em]
& \geq & |V_{2}|+|V_{3}|+|V_{4}|+\frac{2|V_{0}|}{\Delta }+{\frac{1}{\Delta}%
|V_{2}|+|V_{3}|+|V_{4}|} \\[1em]
& = & |V_{2}|+|V_{3}|+|V_{4}|+\frac{2p-2\left(
|V_{2}|+|V_{3}|+|V_{4}|\right) }{\Delta }+{\frac{1}{\Delta}%
|V_{2}|+|V_{3}|+|V_{4}|} \\[1em]
& = & \displaystyle\frac{{2p}+(\Delta -1)\left(
|V_{2}|+|V_{3}|+|V_{4}|\right) }{\Delta }+{\frac{\Delta-1}{\Delta}%
(|V_{3}|+|V_{4}|)} \\[1em]
& \geq & \displaystyle\frac{{2p}+(\Delta -1)\gamma (\Gamma)}{\Delta }%
\end{array}%
\end{equation*}%
\noindent and the result holds. $\ \Box $

\section{{TRD-number} in some families of graphs.}

In this section, we determine exact values for the {TRD-number} in some
particular families of graphs. To do that, we introduce the following
notation: given a positive integer $p\geq 2$ let us denote by $M_{p}$ the
value:

\begin{equation*}
M_{p}=\left\{
\begin{array}{ll}
4\left\lfloor \displaystyle\frac{p}{3}\right\rfloor , & \mbox{if }p\equiv 0
\ (\mathrm{mod \ 3)} \\[1em]
4\left\lfloor \displaystyle\frac{p}{3}\right\rfloor +3, & \mbox{if }p\equiv
1 \ (\mathrm{mod \ 3)} \\[1em]
4\left\lfloor \displaystyle\frac{p}{3}\right\rfloor +4, & \mbox{if }p\equiv
2\ (\mathrm{mod \ 3).}%
\end{array}%
\right.
\end{equation*}

\begin{prop}
\label{paths}Let $p\geq 2$ be a positive integer. Then $\gamma _{\lbrack
3R]}\left( P_{p}\right) =M_{p}.$
\end{prop}

\noindent \textbf{Proof.} Let $P_{p}$ be a path on $p$ vertices with vertex
set $V(P_{p})=\{v_{11},v_{12},v_{13},v_{21},v_{22},$ $v_{23},\ldots
,v_{m1},v_{m2},v_{m3}\}\cup V_{t}^{\prime }$ where $p=3m+t$, with $0\leq
t\leq 2,$ and
\begin{equation*}
V_{t}^{\prime }=\left\{
\begin{array}{ll}
\emptyset & \mbox{if }t=0 \\[0.5em]
\{w_{1}\} & \mbox{if }t=1 \\[0.5em]
\{w_{1},w_{2}\} & \mbox{if }t=2%
\end{array}%
\right.
\end{equation*}%
the set of (ordered) vertices of the path. First of all, let us define a
family of functions $h_{t}$ in $V(P_{p})$ in the following terms,
\begin{equation*}
\begin{array}{c}
h_{t}(v_{ij})=\left\{
\begin{array}{ll}
4 & \mbox{if }j=2 \\[0.5em]
0 & \mbox{if }j\neq 2%
\end{array}%
\right. ,\mbox{for all }i\in \{1,2,\ldots ,m\},t\in \{0,1,2\} \\[2em]
\mbox{and }h_{1}(w_{1})=3,h_{2}(w_{1})=h_{2}(w_{2})=2%
\end{array}%
\end{equation*}%
{Clearly}, $h_{t}$ are 3RDF in $P_{p}$ for all $t=p-3\left\lfloor \frac{p}{3}%
\right\rfloor \in \{0,1,2\}$ and the weight of $h_{t}$ is exactly $M_{p}$.
Therefore, $\gamma _{\lbrack 3R]}\left( P_{p}\right) \leq M_{p}.$

On the other hand, let $h$ be a $\gamma _{\lbrack 3R]}\left( P_{p}\right) $%
-function such that no vertex is {label}ed a 1 under $h$. Let us denote by $%
t=p-3\left\lfloor \frac{p}{3}\right\rfloor \in \{0,1,2\}$ and proceed by
induction in the number of vertices $p$. {Clearly}, $\gamma _{\lbrack
3R]}\left( P_{2}\right) =4\geq M_{2},$ $\gamma _{\lbrack 3R]}\left(
P_{3}\right) =4\geq M_{3}$, $\gamma _{\lbrack 3R]}\left( P_{4}\right) =7\geq
M_{4}$ and $\gamma _{\lbrack 3R]}\left( P_{5}\right) =8\geq M_{5}.$ Let $%
p\geq 6,$ and assume that $\gamma _{\lbrack 3R]}\left( P_{p^{\prime
}}\right) \geq M_{p^{\prime }}$ for all $2\leq p^{\prime }<p.$ Let $%
V(P_{p})=\{y_{1},y_{2},\ldots ,y_{p}\}$ be the set of (ordered) vertices of
the path $P_{p}.$ Consider the following two cases.

\textbf{Case 1. }$h(y_{p-3})\geq 3$ or $h(y_{p-3})\leq 2$ and $h_{|P_{p-3}}$
is a 3RDF.

Then $h_{|P_{p-3}}$ is a 3RDF of $P_{p-3}$ and $%
h(y_{p-2})+h(y_{p-1})+h(y_{p})\geq 4.$ Hence
\begin{equation*}
\gamma _{\lbrack 3R]}\left( P_{p}\right) =h(V(P_{p}))\geq
h_{|P_{p-3}}(V(P_{p-3}))+4\geq M_{p-3}+4=M_{p}.
\end{equation*}

\textbf{Case 2. }$h(y_{p-3})\leq 2$ and $h_{|P_{p-3}}$ is not a 3RDF. {%
Clearly} if $h(y_{p-3})=h(y_{p-2})=0,$ then $h_{|P_{p-3}}$ would be a 3RDF,
which is a contradiction. Thus either $h(y_{p-3})\neq 0$ or $h(y_{p-2})\neq
0.$ Consider the following situations.

\textrm{(a) }$h(y_{p-3})=0$ and $h(y_{p-2})=4.$ Then it is easy to check
that it must be $h(y_{p})+h(y_{p-1})\geq 3$. We define the function $l$ in $%
V(P_{p})$ as follows $l(y_{p})=l(y_{p-2})=0$, $l(y_{p-1})=4$, $l(y_{p-3})=3$
and $l(w)=h(w)$ {\ for the remaining vertices}.

\textrm{(b) }$h(y_{p-3})=0$ and $2\leq h(y_{p-2})\leq 3.$ Then by
considering again that $h(y_{p})+h(y_{p-1})\geq 3$, we may deduce that
\begin{equation*}
\sum_{i=p-4}^{p}h(y_{i})\geq 8.
\end{equation*}%
Define the function $l$ in $V(P_{p})$ as follows $%
l(y_{p})=l(y_{p-2})=l(y_{p-3})=0$, $l(y_{p-1})=l(y_{p-4})=4$ and $l(w)=h(w)$
{\ for the remaining vertices}.

\textrm{(c) }$h(y_{p-3})=2.$ Since $h_{|P_{p-3}}$ is not a 3RDF, we have
that $h(y_{n-2})\geq 2$. Therefore, $\sum_{i=p-3}^{p}h(y_{i})\geq 7.$ We can
define the function $l$ in $V(P_{p})$ as follows $l(y_{p})=l(y_{p-2})=0$, $%
l(y_{p-1})=4$, $l(y_{p-3})=3$ and $l(w)=h(w)$ {\ for the remaining vertices}.

In either case, $l$ is a 3RDF in $P_{p}$ such that $l_{|P_{p-3}}$ is also a
3RDF and hence
\begin{equation*}
\gamma _{\lbrack 3R]}\left( P_{p}\right) =h(V(P_{p}))\geq l(V(P_{p}))\geq
l_{|P_{p-3}}(V(P_{p-3}))+4\geq M_{p-3}+4=M_{p}.\text{ \ }
\end{equation*}%
$\Box $

\begin{prop}
\label{no2chains}Let $\Gamma $ be a connected graph with maximum degree $%
\Delta =2$ and let $h$ be a $\gamma _{\lbrack 3R]}(\Gamma )$-function. If $%
P:y_{1}y_{2}\ldots y_{t}$ is a path in $\Gamma $ such that $h(y_{j})=2$ for
all $j=1,\ldots ,t,$ then $t\leq 3.$ Moreover, there exists a $\gamma
_{\lbrack 3R]}\left( \Gamma \right) $-function such that $h(N(y))=\{0,2\},$
for all vertex $y\in V_{2}$ with degree $2$.
\end{prop}

\noindent \textbf{Proof.} {let $h$ be a $\gamma _{\lbrack 3R]}(\Gamma)$%
-function.} Let $P:y_1y_2\ldots y_t$ be a path in $\Gamma$ {such that} $%
h(y_j)=2$ for all $j=1,\ldots,t. $ If $t\ge 4$ then the function $l$ defined
as follows: $l(y_1)=3, l(y_2)=0$ and $l(y^{\prime})=h(y^{\prime}) $ {for
each $y^{\prime}\in V(\Gamma)-\{y_1,y_2\}$} would be a 3RDF in $\Gamma$ with
{$l(V(\Gamma))={h(V(\Gamma))}-1$}, which is impossible. So, $%
t\le3.$

On the other hand, let $l$ be a $\gamma _{\lbrack 3R]}(\Gamma )$-function
that minimizes the size of the set $V_{2}$ of vertices assigned a 2 under $h$%
. {Let }$y\in V_{2}$ {be a vertex of degree two and let }$N(y)=\{v,w\}$. {If
}$N(y)\subseteq V_{2},$ then the function $l$ defined in this way: {$l(y)=0$%
, $l(v)=l(w)=3$ and $l(z)=h(z)$ for each vertex $z\in V(\Gamma )-\{y,v,w\}$}
would be a $\gamma _{\lbrack 3R]}(\Gamma )$-function with fewer vertices
assigned {a $2$ than $h.$ Hence we can assume, without loss of generality,
that }$v\in V_{0}$ and $w\in V_{3}\cup V_{4}.$ {By the definition of a 3RDF,
}$v$ must have either one neighbour in $V_{4}$, or either one neighbour in $%
V_{2}$ and another one in $V_{3}$ or either three neighbours in $V_{2}.$ If $%
N(v)\cap V_{4}\neq \emptyset ,$ {then reassigning }$y$ the value $1$ instead
of $2$ provides a 3RDF with weight $h(\Gamma )-1,$ a contradiction. If $v$
has one neighbour in $V_{2}$ and another one $z$ in $V_{3},$ {then
reassigning vertices }$y,v,w,z$ the values $0,0,4,4, $respectively provides
a $\gamma _{\lbrack 3R]}(\Gamma ) $-function with fewer vertices assigned {a
$2$ than $h.$ Finally, assume that }$v$ has three neighbours in $V_{2}.$ {%
Then reassigning vertices }$y,v,w$ the values $0,1,4,$ respectively provides
a $\gamma _{\lbrack 3R]}(\Gamma )$-function with fewer vertices assigned {a $%
2$ than $h.$ } 

\begin{prop}
\label{cycles}For any integer {$p\geq 3,$}
\begin{equation*}
\gamma _{\lbrack 3R]}\left( C_{p}\right) =\left\{
\begin{array}{ll}
\displaystyle\left\lceil \frac{4p}{3}\right\rceil & \mbox{if either }%
p=4,5,7,10\mbox{ or }{p\equiv 0\ (\mathrm{mod\ 3).}} \\[2em]
\displaystyle\left\lceil \frac{4p}{3}\right\rceil +1 & \mbox{if }p\neq
4,5,7,10\mbox{ and }p\equiv 1,2\ {(\mathrm{mod\ 3)}}.%
\end{array}%
\right.
\end{equation*}
\end{prop}

\noindent \textbf{Proof.} First of all, observe that $M_{p}=\left\lceil
\frac{4p}{3}\right\rceil +t$ where $t=0,$ if $p\equiv 0(\mod 3);$ and $t=1,$
if $p>0(\mod 3).$

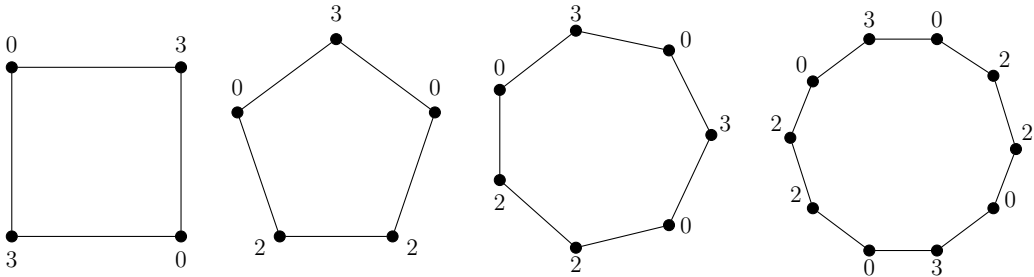
\begin{figure}[h]
\centering
\begin{tikzpicture}[scale=.75, transform shape]


\node [draw, shape=circle,fill=black,scale=0.5] (u0) at  (0,2) {};
\node [draw, shape=circle,fill=black,scale=0.5] (u1) at  (0,5) {};
\node [draw, shape=circle,fill=black,scale=0.5] (u2) at  (3,5) {};
\node [draw, shape=circle,fill=black,scale=0.5] (u3) at  (3,2) {};
\draw (u0)--(u1);
\draw (u1)--(u2);
\draw (u2)--(u3);
\draw (u3)--(u0);
\node[draw=none, fill=none] at (0,1.6){ $3$};
\node[draw=none, fill=none] at (0,5.4){ $0$};
\node[draw=none, fill=none] at (3,5.4){ $3$};
\node[draw=none, fill=none] at (3,1.6){ $0$};

\node [draw, shape=circle, fill=black, scale=0.5] (a1) at  (5.75 ,5.5 ) {};
\node [draw, shape=circle, fill=black, scale=0.5] (a2) at  (4 ,4.2) {};
\node [draw, shape=circle, fill=black, scale=0.5] (a3) at  (7.5 ,4.2) {};
\node [draw, shape=circle, fill=black, scale=0.5] (a4) at  (4.75 ,2 ) {};
\node [draw, shape=circle, fill=black, scale=0.5] (a5) at  (6.75 ,2 ) {};
\node [draw=none, fill=none] at  (5.75 , 5.95 ) {3};
\node [draw=none, fill=none] at  (4 ,4.65) {0};
\node [draw=none, fill=none] at  (7.5 ,4.65 ) {0};
\node [draw=none, fill=none] at  (4.4 ,1.8 ) {2};
\node [draw=none, fill=none] at  (7.1,1.8) {2};
\draw (a1)--(a2)--(a4)--(a5)--(a3)--(a1);

\node [draw, shape=circle, fill=black, scale=0.5] (b1) at  (10 ,5.65 ) {};
\node [draw, shape=circle, fill=black, scale=0.5] (b2) at  (8.65 ,4.6 ) {};
\node [draw, shape=circle, fill=black, scale=0.5] (b3) at  (8.65 ,3 ) {};
\node [draw, shape=circle, fill=black, scale=0.5] (b4) at  (10 ,1.8 ) {};
\node [draw, shape=circle, fill=black, scale=0.5] (b5) at  (11.65 ,2.2) {};
\node [draw, shape=circle, fill=black, scale=0.5] (b6) at  (12.4 ,3.8 ) {};
\node [draw, shape=circle, fill=black, scale=0.5] (b7) at  (11.65 ,5.3 ) {};
\node [draw=none, fill=none] at  (10 , 5.95 ) {3};
\node [draw=none, fill=none] at  (8.65 , 5 ) {0};
\node [draw=none, fill=none] at  (8.65 , 2.6 ) {2};
\node [draw=none, fill=none] at  (10 , 1.5 ) {2};
\node [draw=none, fill=none] at  (11.95 , 2.2 ) {0};
\node [draw=none, fill=none] at  (11.95 , 5.5 ) {0};
\node [draw=none, fill=none] at  (12.65 , 4 ) {3};

\draw (b1)--(b2)--(b3)--(b4)--(b5)--(b6)--(b7)--(b1);

\node [draw, shape=circle, fill=black, scale=0.5] (c1) at  (15.2 , 1.75 ) {};
\node [draw, shape=circle, fill=black, scale=0.5] (c2) at  (16.4 , 1.75 ) {};
\node [draw, shape=circle, fill=black, scale=0.5] (c3) at  (17.4 , 2.5 ) {};
\node [draw, shape=circle, fill=black, scale=0.5] (c4) at  (17.8 , 3.55 ) {};
\node [draw, shape=circle, fill=black, scale=0.5] (c5) at  (17.4 , 4.85 ) {};
\node [draw, shape=circle, fill=black, scale=0.5] (c6) at  (16.4 , 5.5 ) {};
\node [draw, shape=circle, fill=black, scale=0.5] (c7) at  (15.2 , 5.5 ) {};
\node [draw, shape=circle, fill=black, scale=0.5] (c8) at  (14.2 , 4.75 ) {};
\node [draw, shape=circle, fill=black, scale=0.5] (c9) at  (13.8 , 3.75 ) {};
\node [draw, shape=circle, fill=black, scale=0.5] (c10) at  (14.2 , 2.5 ) {};

\node [draw=none, fill=none] at  (15.2 , 1.45 ) {0};
\node [draw=none, fill=none] at  (16.4 , 1.45) {3};
\node [draw=none, fill=none] at  (17.7, 2.65 ) {0};
\node [draw=none, fill=none] at  (18 , 3.85 ) {2};
\node [draw=none, fill=none] at  (17.6, 5.15 ) {2};
\node [draw=none, fill=none] at  (16.4 , 5.85 ) {0};
\node [draw=none, fill=none] at  (15.2 , 5.85 ) {3};
\node [draw=none, fill=none] at  (14 , 5.05 ) {0};
\node [draw=none, fill=none] at  (13.55 , 4 ) {2};
\node [draw=none, fill=none] at  (13.9 , 2.75 ) {2};

\draw (c1)--(c2)--(c3)--(c4)--(c5)--(c6)--(c7)--(c8)--(c9)--(c10)--(c1);

\end{tikzpicture}
\caption{{3RDF}{\protect\small s for $C_{4},C_{5},C_{7},C_{10}$} }
\label{partcycles}
\end{figure}

Let us begin with the cycles $C_{4},C_{5},C_{7}$ and $C_{10}$ by considering
the {3RDF} depicted in Figure \ref{partcycles}.

Hence, we deduce that $\gamma _{\lbrack 3R]}(C_{4})\leq 6,\gamma _{\lbrack
3R]}(C_{5})\leq 7,\gamma _{\lbrack 3R]}(C_{7})\leq 10$ and $\gamma _{\lbrack
3R]}(C_{10})\leq 14.$ Now, let $C_{p}:y_{1}y_{2}\ldots y_{p}y_{1}$ be a
cycle, with $p\neq 4,5,7,10.$ Since each $\gamma _{\lbrack 3R]}$-function in
a path induces a 3RDF in the corresponding cycle, then $\gamma _{\lbrack
3R]}(C_{p})\leq \gamma _{\lbrack 3R]}(C_{p}-y_{1}y_{2})=M_{p}$, by applying
the Proposition \ref{paths}.

On the other hand, let $h=(V_0,\emptyset,V_2,V_3,V_4)$ be a $%
\gamma_{[3R]}(C_p)$-function. If there exists a vertex $y_i$ such that
either $h(y_i)=h(y_{i+1})=0$ or either $h(y_{i+1})=0, h(y_{i+2})=4$ or
either $h(y_i),h(y_{i+1})\ge 3$ (where the subscripts are considered modulus
$p$), then $h_{|C_p-y_iy_{i+1}}$ is a 3RDF in the path $C_p-y_iy_{i+1}$.
Therefore,
\begin{equation}  \label{direct}
\gamma_{[3R]}(C_p)={h(V(C_p))}=w(h_{|C_p-y_iy_{i+1}})\ge
\gamma_{[3R]}(C_p-y_iy_{i+1})=M_p.
\end{equation}

So, {assume} that $h=(V_{0},\emptyset ,V_{2},V_{3},\emptyset )$ is a $\gamma
_{\lbrack 3R]}(C_{p})$-function such that $E[V_{0},V_{0}]=E[V_{3},V_{3}]=%
\emptyset ,$ that is sets $V_{0}$ and $V_{3}$ are independent. Let $r$
(resp. $s,t$) be the number of vertices in $C_{p}$ {label}ed with the label $%
3$ (resp. $2,0$) {under }$h$. {Clearly}, $r+s+t=p$, $3r+2s={%
h(V(C_p))}$ and, since each $0$ must be {connecting to} {at least} a $3,$ we
have that $t\leq 2r.$

Besides, as ${h(V(C_p))}$ is minimum, if $h(y)=3,$ then either $%
h(N(y))=\{0\}$ or either $h(N(y))=\{0,2\}.$ Next, we shall show that it is
possible to have $E[V_{2},V_{3}]=\emptyset .$ {Assume} that there are two
adjacent vertices $y_{i},y_{i+1}$ such that $h(y_{i})=2$ and $h(y_{i+1})=3.$
In this situation, it must be $h(y_{i+2})=0,$ (where all the subscripts are
considered to be modulus $p$). It is straightforward to see that $%
h(y_{i-1})=0,$ because {\ for the remaining vertices} $%
h(y_{i-1})+h(y_{i})+h(y_{i+1})\geq 7$ and the function $l(y_{i-1})=3,$ $%
l(y_{i})=0$ and {$l(y^{\prime})=h(y^{\prime})$} for all $y^{\prime}\neq
y_{i-1},y_{i}$ would be a 3RDF with weight less than ${h(V(C_p))}%
,$ which is impossible. Since $h(y_{i-1})=0$ then it must be $h(y_{i-2})=3$
because $h(y_{i})=2,$ which implies that $h(y_{i-3})\leq 2$ due to $%
E[V_{3},V_{3}]=\emptyset .$ Consider the two possible situations.

If $h(y_{i-3})=2$ then, reasoning analogously to when we showed that it must
necessarily be $h(y_{i-1})=0,$ we have that $h(y_{i-4})=0.$ Hence, the
function $l(y_{i-3})=4$, $l(y_{i-2})=0$, $l(y_{i-1})=3$, $l(y_{i})=0$ and $%
l(y^{\prime})=h(y^{\prime})$ {\ for the remaining vertices} is a 3RDF with ${%
h(V(C_p))=l(V(C_p))}.$ Besides, since $l(y_{i-3})=4 $ and $%
l(y_{i-2})=0$, by applying $(\ref{direct}),$ we have that ${%
h(V(C_p))=l(V(C_p))}\geq M_{p}.$

If $h(y_{i-3})=0$ then $h(y_{i-4})\in \{2,3\}.$ In case that $h(y_{i-4})=3$
we consider the function $l(y_{i-4})=4$, $l(y_{i-3})=0$,$l(y_{i-2})=0$, $%
l(y_{i-1})=4$, $l(y_{i})=0$ and $l(y^{\prime})=h(y^{\prime})$ {\ for the
remaining vertices}. {Clearly}, $l$ is a 3RDF with ${%
h(V(C_p))=l(V(C_p))}$ and since $l(y_{i-2})=0,l(y_{i-1})=4$ we may apply $(%
\ref{direct})$ to deduce that ${h(V(C_p))=l(V(C_p))}\geq M_{p}.$
Hence let $h(l_{i-4})=2. $ Then either $h(y_{i-5})=2$ and $h(y_{i-6})=0$ or
either $h(y_{i-5})=3$ and $h(y_{i-6})=0$. In the former situation, let us
consider the function $l(y_{i-5})=3$, $l(y_{i-4})=0$,$l(y_{i-3})=3$, $%
l(y_{i-2})=0$, $l(y_{i-1})=3$, $l(y_{i})=0$ and $l(y^{\prime})=h(y^{\prime})$
{\ for the remaining vertices}. Again it is easy to deduce that ${%
h(V(C_p))=l(V(C_p))}\geq M_{p}.$ In the latter situation, we
define the function $l$ as follows: $l(y_{i-4})=0$,$l(y_{i-3})=3$, $%
l(y_{i-2})=0$, $l(y_{i-1})=4$, $l(y_{i})=0$ and $l(y^{\prime})=h(y^{\prime})$
{\ for the remaining vertices}. So, $l$ is a 3RDF with ${%
h(V(C_p))=l(V(C_p))}$ and the joining of a vertex with label $2$ with a
vertex {label}ed with a $3$ is avoided.

Summing up, without loss of generality, we may suppose that $%
E(G)=E[V_{0},V_{2}]\cup E[V_{2},V_{2}]\cup E[V_{0},V_{3}].$ Moreover, since
each vertex {label}ed with a $3$ is {connecting to} two vertices labeled
with $0, $ then $|E[V_{0},V_{3}]|=2r.$ According to Proposition \ref%
{no2chains}, it is clear that $h(N(y))=\{0,2\}$ for all vertex $y\in V_{2}$.
Therefore we derive that $E[V_{2},V_{2}]=|V_{2}|/2=s/2,$ and $%
|E[V_{0},V_{2}]|=|V_{2}|=s.$ The latter lead us to the following relations
between $r,s,t$
\begin{equation*}
\left\{
\begin{array}{rcrcrcl}
r & + & s & + & t & = & p \\
4r & + & 3s &  &  & = & 2p \\
3r & + & 2s &  &  & = & w(f)%
\end{array}%
\right.
\end{equation*}%
from which we deduce that: $r=3{h(V(C_p))}-4p,s=6p-4{%
h(V(C_p))},t={h(V(C_p))}-p.$ Let us note that in this situation $%
r=0$ implies $t=0,$ because $t\leq 2r$ and hence, by Proposition \ref%
{no2chains}, we may assume that $r>0.$ If $p=3m,$ then we have that $\gamma
_{\lbrack 3R]}(C_{p})=w(f)>4p/3=M_{p},$ a contradiction. If $p=3m+1$ then it
must be ${h(V(C_p))}>\frac{4}{3}(3m+1)=4m+4/3$ and therefore ${%
h(V(C_p))}\geq 4m+2$. But, if ${h(V(C_p))}=4m+2$
then $r=2,s=2m-2$ and $t=m+1\leq 2r=4$ which implies that $1\leq m\leq 3.$
For $m=1$ ($m=2,3$ respectively) we obtain that $\gamma _{\lbrack
3R]}(C_{4})=6,(\gamma _{\lbrack 3R]}(C_{7})=10,\gamma _{\lbrack
3R]}(C_{10})=14$, respectively) as we needed to show. Therefore, $\gamma
_{\lbrack 3R]}(C_{p})={h(V(C_p))}=4m+3=\left\lceil \frac{4p}{3}%
\right\rceil +1=M_{p},$ for all $p=3m+1$ with $p\neq 4,7,10$.

If $p=3m+2$ then it must be ${h(V(C_p))}>\frac{4}{3}(3m+2)=4m+8/3
$ and therefore $w(h)\geq 4m+3$. Now, if ${h(V(C_p))}=4m+3$ then
$r=1,s=2m$ and $t=m+1\leq 2r=2$ which implies that $m=1,$ and hence, we
obtain that $\gamma _{\lbrack 3R]}(C_{5})=7$. Besides, we deduce that $%
\gamma _{\lbrack 3R]}(C_{p})={h(V(C_p))}=4m+4=\left\lceil \frac{%
4p}{3}\right\rceil +1=M_{p},$ for all $p=3m+2$ with $p\neq 5.$ $\ \Box $

Finally, we give two upper bounds of $\gamma _{\lbrack 3R]}\left(
\Gamma\right) $ in terms of the diameter and the girth of the graphs that
are immediate consequences of Propositions \ref{paths} and \ref{cycles}.

\begin{prop}
\label{diam}Let $\Gamma $ be {an ntc-graph} with $p$ vertices and diameter $%
\mathrm{diam(\Gamma )}$. Then
\begin{equation*}
\gamma _{\lbrack 3R]}\left( \Gamma \right) \leq 3p-\frac{5\mathrm{diam}%
(\Gamma )}{3}+\frac{7}{3}
\end{equation*}
\end{prop}

{\ \noindent \textbf{Proof.} Let us denote by $d=\mathrm{diam}(\Gamma)$ and
let $P_{d+1}$ be a diametral path with $d+1$ vertices. By Proposition \ref%
{paths}, it is easy to check that $\gamma_{[3R]}\left( P_{d+1} \right) \le
4\left\lfloor \frac{d+1}{3}\right\rfloor +4.$ Any $\gamma_{[3R]}\left(
P_{d+1} \right)$-function can be extended to a 3DRF in $\Gamma$ by assigning
a $3$ to every vertex in $V(\Gamma) \setminus V(P_{d+1}).$ Hence
\begin{equation*}
\begin{array}{rcl}
\gamma_{[3R]}\left( \Gamma \right) & \le & \gamma_{[3R]}\left( P_{d+1}
\right)+ 3(p-d-1) \\[1em]
& \le & 4\left\lfloor \displaystyle\frac{d+1}{3}\right\rfloor +4 +3p-3d-3 \\%
[1em]
& \le & \displaystyle\frac{4}{3} d+\frac{4}{3}+4+3p-3d-3 = 3p-\frac{5 d}{3}+%
\frac{7}{3}%
\end{array}%
\end{equation*}
$\hfill \Box$ }

\begin{prop}
\label{girth}Let $\Gamma $ be {an ntc-graph} of order $p$ and with girth $%
g(\Gamma )$. Then
\begin{equation*}
\gamma _{\lbrack 3R]}\left( \Gamma \right) \leq 3p+2-\frac{5g(\Gamma )}{3}.
\end{equation*}
\end{prop}

{\ \noindent \textbf{Proof.} Let us denote by $g=\mathrm{girth}(\Gamma )$
and let $C_{g}$ be a cycle in $\Gamma $ with $g$ vertices. By Proposition %
\ref{cycles}, we have that $\gamma _{\lbrack 3R]}\left( C_{g}\right) \leq
\left\lceil \frac{4g}{3}\right\rceil +1.$ Any $\gamma _{\lbrack 3R]}\left(
C_{g}\right) $-function can be extended to a 3DRF in $\Gamma $ by assigning
a $3$ to every vertex in $V(\Gamma )\setminus V(C_{g}). $ Hence
\begin{equation*}
\begin{array}{rcl}
\gamma _{\lbrack 3R]}\left( \Gamma \right) & \leq & \gamma _{\lbrack
3R]}\left( C_{g}\right) +3(p-g)\leq \left\lceil \displaystyle\frac{4g}{3}%
\right\rceil +1+3p-3g \\[1em]
& \leq & \displaystyle\frac{4g}{3}+2+3p-3g=3p-\frac{5g}{3}+2.%
\end{array}%
\end{equation*}%
$\hfill \Box $ }

\textbf{Acknowledgements}\newline

The authors are grateful to anonymous referee for his/her remarks and
suggestions that helped improve the manuscript.
H. Abdollahzadeh Ahangar was supported by the Babol Noshirvani University of
Technology under research grant number BNUT/385001/99.


\end{document}